\documentclass[12pt]{amsart}
\usepackage{amsbsy}
\usepackage{graphicx,epsfig,subfigure,psfrag}
\begin{document}

\newtheorem{theorem}{Theorem}
\newtheorem{proposition}{Proposition}
\newtheorem{lemma}{Lemma}
\newtheorem{corollary}{Corollary}
\newtheorem{definition}{Definition}
\newtheorem{remark}{Remark}
\newcommand{\be}{\begin{equation}}
\newcommand{\ee}{\end{equation}}
\newcommand{\tex}{\textstyle}
\numberwithin{equation}{section} \numberwithin{theorem}{section}
\numberwithin{proposition}{section} \numberwithin{lemma}{section}
\numberwithin{corollary}{section}
\numberwithin{definition}{section} \numberwithin{remark}{section}
\newcommand{\ren}{\mathbb{R}^N}
\newcommand{\re}{\mathbb{R}}
\newcommand{\n}{\nabla}
\newcommand{\iy}{\infty}
\newcommand{\pa}{\partial}
\newcommand{\fp}{\noindent}
\newcommand{\ms}{\medskip\vskip-.1cm}
\newcommand{\mpb}{\medskip}
\newcommand{\BB}{{\bf B}}
\newcommand{\AAA}{{\bf A}}
\newcommand{\Am}{{\bf A}_{2m}}
\newcommand{\ef}{\eqref}
\newcommand{\eee}{{\mathrm e}}
\newcommand{\ii}{{\mathrm i}}
\renewcommand{\a}{\alpha}
\renewcommand{\b}{\beta}
\newcommand{\g}{\gamma}
\newcommand{\G}{\Gamma}
\renewcommand{\d}{\delta}
\newcommand{\D}{\Delta}
\newcommand{\e}{\varepsilon}
\newcommand{\var}{\varphi}
\renewcommand{\l}{\lambda}
\renewcommand{\o}{\omega}
\renewcommand{\O}{\Omega}
\newcommand{\s}{\sigma}
\renewcommand{\t}{\tau}
\renewcommand{\th}{\theta}
\newcommand{\z}{\zeta}
\newcommand{\wx}{\widetilde x}
\newcommand{\wt}{\widetilde t}
\newcommand{\noi}{\noindent}
\newcommand{\uu}{{\bf u}}
\newcommand{\UU}{{\bf U}}
\newcommand{\VV}{{\bf V}}
\newcommand{\ww}{{\bf w}}
\newcommand{\vv}{{\bf v}}
\newcommand{\WW}{{\bf W}}
\newcommand{\hh}{{\bf h}}
\newcommand{\di}{{\rm div}\,}
\newcommand{\inA}{\quad \mbox{in} \quad \ren \times \re_+}
\newcommand{\inB}{\quad \mbox{in} \quad}
\newcommand{\inC}{\quad \mbox{in} \quad \re \times \re_+}
\newcommand{\inD}{\quad \mbox{in} \quad \re}
\newcommand{\forA}{\quad \mbox{for} \quad}
\newcommand{\whereA}{,\quad \mbox{where} \quad}
\newcommand{\asA}{\quad \mbox{as} \quad}
\newcommand{\andA}{\quad \mbox{and} \quad}
\newcommand{\withA}{,\quad \mbox{with} \quad}
\newcommand{\orA}{,\quad \mbox{or} \quad}
\newcommand{\ssk}{\smallskip}
\newcommand{\LongA}{\quad \Longrightarrow \quad}
\def\com#1{\fbox{\parbox{6in}{\texttt{#1}}}}
\def\N{{\mathbb N}}
\def\A{{\cal A}}
\def\WW{{\cal W}}
\newcommand{\de}{\,d}
\newcommand{\eps}{\varepsilon}
\newcommand{\spt}{{\mbox spt}}
\newcommand{\ind}{{\mbox ind}}
\newcommand{\supp}{{\mbox supp}}
\newcommand{\dip}{\displaystyle}
\newcommand{\prt}{\partial}
\renewcommand{\theequation}{\thesection.\arabic{equation}}
\renewcommand{\baselinestretch}{1.2}

\title
{\bf Towards the KPP--problem and
  ${\bf {log\, t}}$--front \\ shift for higher-order nonlinear  PDEs III.
\\
 Dispersion and hyperbolic equations}

\author{
V.A.~Galaktionov}


\address{Department of Mathematical Sciences, University of Bath,
 Bath BA2 7AY, UK}
\email{masvg@bath.ac.uk}



  \keywords{KPP-problem, travelling wave, stability,  higher-order
  PDEs, odd-order dispersion and hyperbolic equations, front $\log t$-shift}
 \subjclass{35K55, 35K40, 35K65}
 \date{\today}




\begin{abstract}

Some aspects of extensions of ideas of Kolmogorov, Petrovskii, and
Piskunov (1937) \cite{KPP} on travelling wave propagation in the
reaction-diffusion equation
 $$
 u_t=u_{xx}+u(1-u) \inB \re \times \re_+, \quad u_0(x)=H(-x) \equiv\{1\,\,
 \mbox{for}\,\, x<0; \,\,\, 0 \,\, \mbox{for} \,\, x \ge 0\},
 $$
are discussed. The
 present paper continues the study began in \cite{GKPPI, GKPPII}
 for higher-order
 parabolic
 {\em semilinear} and {\em quasilinear} {\em bi-harmonic
 equations} such as
 $$
 u_t= -u_{xxxx} +u(1-u), \quad
u_t=-(|u|^n u)_{xxxx}+u(1-u) \quad (n>0),
 \quad \mbox{etc.}
 $$

Here, higher-order {\em dispersion} equations such as ($D_x=
\frac{\partial}{\partial x}$ and $D_t= \frac{\partial}{\partial
t}$)
 $$
   u_{t}=- D_x^{11}u +u(1-u) \quad \mbox{and up to} \quad
 D_t^9 u =- D_x^{11} u + u(1-u), \quad \mbox{etc.}
 $$
 are studied. Some features of KPP-like
 results are also shown to exist  for semilinear {\em dispersion-parabolic} equations
 such as
 $$
 u_{ttt}= - D_x^{10}u +u(1-u), \quad u_{ttttt}=D_x^{10}u+u(1-u), \quad \mbox{and others},
  $$
  and for pure {\em hyperbolic} ones
 $$
 u_{tt}=-u_{xxxx}+u(1-u) \quad \mbox{and up to} \quad u_{tttt}=-D_x^{10} u +u(1-u), \quad
 \mbox{etc.}
 $$
 As an example, we also treat a {\em quasilinear} PDE
 $u_t=-D_x^{11}(|u|^n u)+u(1-u)$, with $n>0$.
 Two main questions  are:
 (i) existence of travelling waves via any analytical/numerical methods, and
 (ii) their stability and
 derivation of the
 $\log t$-shifting of moving fronts.

\end{abstract}

\maketitle



\setcounter{equation}{0}
\section{Introduction: the classic KPP-problem and other higher-order PDE models}
 \label{Sect1}
  \setcounter{equation}{0}











\subsection{The classic KPP-problem of 1937}

The classic KPP-problem \cite{KPP} (1937)
 \be
 \label{1.1}
  u_t = u_{xx} +u(1-u) \inB \re \times \re_+,
\quad u(x,0)=u_0(x) \,\,\, \mbox{in \,\,$\re$},
 \ee
    with
the step (Heaviside) initial function
 \be
 \label{1.H}
u_0(x) =H(-x) \equiv \left\{ \begin{matrix}1, \quad x<0;\\ 0,
\quad x\ge 0,
\end{matrix}
\right.
 \ee
consists of studying large-time convergence of the solution of the
Cauchy problem \ef{1.1}, \ef{1.H} to the unique {\em minimal} {\em
travelling wave} solution, with the minimal speed $\l_0=2$ of the
standard form
 \be
 \label{TW1}
 \left\{
  \begin{matrix}
 u_*(x,t)=f(y), \,\, y=x-\l_0 t, \,\,\,\,\mbox{where}\qquad\qquad\qquad\qquad\quad\ssk\ssk\\
  -\l_0 f'=f''+f(1-f), \,\, y \in \re; \quad
 f(-\iy)=1, \,\, f(+\iy)=0.
 \end{matrix}
 \right.
  \ee
The KPP paper \cite{KPP} contains a number of pioneering
remarkable results, which founded several new directions of modern
nonlinear PDE theory. For further use, it suffices for us to refer
to a survey and key references in \cite{GKPPI}.

Let us state the main result of \cite{KPP}.
The convergence to the minimal TW \ef{TW1} was performed in the TW
moving frame, proving that
  the
TW front moves like
 \be
 \label{1.2}
  x_f(t) = 2t -
g(t) \quad \mbox{as $t\to \infty$}, \quad \mbox{with $g(t)=o(t)$},
 \ee
 where the front location
 $x_f(t)$ is
uniquely determined from the equation
 \be
 \label{TW2}
  \tex{
 u(x_f(t),t) = \frac 12 \quad \mbox{for all $t \ge 0$}.
 }
  \ee
Then the convergence result of \cite{KPP} takes the form:
 \be
 \label{TW3}
 u(x_f(t)+y,t) \to f(y) \asA t \to + \iy \quad \mbox{uniformly in
 $y \in \re$}.
 \ee

In 1983,
 Bramson \cite{Br},
 using
probabilistic techniques,
proved that  there exists {\em unbounded} $\log t$-shift of the
moving TW front
 \be
 \label{1.3}
  g(t) = k \log t(1+o(1)), \quad
\mbox{with $ k= {\frac 32}$},
  \ee
  Therefore, \ef{1.3} implies eventual, as $t \to +\iy$, {\em infinite} retarding of the
  solution $u(x,t)$ from the corresponding minimal TW, thought the
  convergence \ef{TW3} takes place in the TW frame.

\subsection{KPP-like problem
to higher-order semilinear and quasilinear  parabolic PDEs}

We dealt with such semilinear reaction-diffusion PDEs in
\cite{GKPPI}, where the main basic model was
 the {\em semilinear bi-harmonic equation}, i.e., a fourth-order semilinear heat equation (SHE--4)
 \be
 \label{E4}
 u_t= -u_{xxxx} +u(1-u) \inB \re \times \re_+.
  \ee
The corresponding TW with the speed of propagation $\l$ is then
governed by the following fourth-order ODE:
 \be
 \label{E5}
u_*(x,t)=f(y), \quad  y=x - \l t \LongA -\l f'=-f^{(4)}+f(1-f),
 \ee
 with the singular boundary conditions at infinity:
 \be
 \label{BC1}
 f(y) \to 0  \andA f(y) \to 1 \asA y \to \pm \iy \quad
 \mbox{``maximally" exponentially fast.}
  \ee
 In particular, we have found a ``maximal speed $\l_{\rm max}=1.27148...$ such that
 \be
 \label{SpeedMax}
 \mbox{TW profiles $f(y;\l)$ exist for all $0< \l < \l_{\rm max}$, and nonexistent
 for $\l>\l_{\rm max}$}.
 \ee

In \cite{GKPPII}, we extend some of the above results to the {\em
quasilinear} KPP--$4n$ problem for
\be
 \label{E4n}
 u_t= -(|u|^n u)_{xxxx} +u(1-u) \inB \re \times \re_+,
  \ee
  where $n>0$ is a parameter, as well as to some other parabolic
  equations.

 In a similar manner, we studied in \cite{GKPPI} the {\em
semilinear tri-harmonic equation} (SHE-6):
 \be
 \label{E6}
 u_t= u_{xxxxxx} +u(1-u) \inB \re \times \re_+, \quad u(x,0)=H(-x)
 \inB \re.
  \ee
The corresponding TW with the speed of propagation $\l$ is then
governed by the following fourth-order ODE:
 \be
 \label{E7}
u_*(x,t)=f(y), \quad  y=x - \l t \LongA -\l f'=f^{(6)}+f(1-f),
\,\,\, \mbox{with (\ref{BC1})}.
 \ee
The \ef{SpeedMax} remains true, with $\l_{\rm max}=2.12110...$\,.

In general, in \cite{GKPPI}, we presented some numerical evidence
on existence of various TWs and there properties for semilinear
parabolic $2m$th-order PDEs (SHE-$2m$) such as (here, $D_x=
\frac{\partial}{\partial x}$)
\be
\label{m1}
 \tex{
 u_t= (-1)^{m+1} D_x^{2m} u + u(1-u) \quad \big(D_x= \frac {\mathrm d}{{\mathrm d}x}\big), \quad \mbox{with the ODEs}
 }
 \ee
 \be
 \label{m2}
u_*(x,t)=f(y), \quad  y=x - \l t \LongA -\l
f'=(-1)^{m+1}f^{(2m)}+f(1-f)
 \ee
 (plus
(\ref{BC1})), and rather sharply estimated $\l_{\rm max}=\l_{\rm
max}(m)>0$ for $m=3,4,5$,
 i.e., up to the tenth-order parabolic equation as in  \ef{m1}.


 \subsection{Results I: dispersion  PDEs (Section \ref{SDisp1})}

In the present paper, we will deal with higher-order {\em
dispersion}\footnote{Here, we use a PDE classification, associated
with some {\em a priori} bounds admitted by the principal linear
differential operators by multiplying by appropriate time
derivatives $D^k_tu$ in the $L^2$-metric. Therefore, while, for
most of parabolic and hyperbolic equations, such a classification
coincides with the classic rigorous Petrovskii's one \cite{Pet61}
(recall: ``parabolic in Petrovskii's sense"), for others, our
classification can be different and uses a ``more applied"
understanding of dispersion and related phenomena. Anyway, we do
not think that, for truly higher-order (in both $x$ and $t$ up to
11th or 12th orders) semilinear PDEs, any formal classification
may somehow essentially help to understand the nature of TW
patterns obtained below (note that such ``evolution" patterns have
been obtained even for obviously {\em elliptic} PDEs, for which a
standard evolution interpretation makes no sense, due to
Hadamard's example of an ill-posed Cauchy problem).}
 equations
 \be
 \label{m31}
   u_{t}= - D_x^{11}u +u(1-u) \LongA - \l f'= -f^{(11)} + f(1-f);
  \ee
\be
 \label{m33}
   u_{ttt}= - D_x^{11}u +u(1-u) \LongA - \l^3 f'''= -f^{(11)} + f(1-f);
  \ee
  \be
 \label{m35}
   u_{ttttt}= - D_x^{11}u +u(1-u) \LongA - \l^5 f^{(5)}= -f^{(11)} + f(1-f);
  \ee
  \be
 \label{m37}
   D_t^7 u = - D_x^{11}u +u(1-u) \LongA - \l^7 f^{(7)}= -f^{(11)} + f(1-f);
  \ee
\be
 \label{m39}
   D_t^9 u = - D_x^{11}u +u(1-u) \LongA - \l^9 f^{(9)}= -f^{(11)} + f(1-f);
  \ee
plus singular boundary conditions \ef{BC1}. We have chosen rather
higher-order dispersion operator $D_x^{11}$ on the right-hand side
in order to avoid a kind of ``temptation" to rely on linearized
analysis, which was heavily used in \cite[\S~2]{GKPPI} in studying
the lower-order bi-harmonic KPP-problem for \ef{E4}.

 \subsection{Results II: dispersion-hyperbolic  PDEs (Section \ref{S.D-H})}

We next consider existence of TW solutions of higher-order {\em
dispersion-hyperbolic} equations
 \be
 \label{m32dh}
   u_{tt}= - D_x^{11}u +u(1-u) \LongA  \l^2 f''= -f^{(11)} + f(1-f);
  \ee
\be
 \label{m34dh}
   u_{tttt}= - D_x^{11}u +u(1-u) \LongA  \l^4 f^{(4)}= -f^{(11)} + f(1-f);
  \ee
  \be
 \label{m36dh}
   u_{tttttt}= - D_x^{11}u +u(1-u) \LongA  \l^6 f^{(6)}= -f^{(11)} + f(1-f);
  \ee
  \be
 \label{m38dh}
   D_t^8 u = - D_x^{11}u +u(1-u) \LongA  \l^8 f^{(8)}= -f^{(11)} + f(1-f);
  \ee
\be
 \label{m310dh}
   D_t^{10} u = - D_x^{11}u +u(1-u) \LongA  \l^{10} f^{(10)}= -f^{(11)} + f(1-f);
  \ee
plus singular boundary conditions \ef{BC1}.

 \subsection{Results III: dispersion-parabolic  PDEs (Section \ref{S.D-P})}

We next consider, in the KPP setting, the following equations,
which, in view of certain estimates, can be considered as  {\em
dispersion-parabolic} (the signs in front of $D_x^{10}$ avoids
``backward parabolic" features) models
 \be
 \label{m33par}
   u_{ttt}= - D_x^{10}u +u(1-u) \LongA - \l^3 f'''=- f^{(10)} + f(1-f);
  \ee
\be
 \label{m35par}
   u_{ttttt}=  D_x^{10}u +u(1-u) \LongA - \l^5 f^{(5)}= f^{(10)} + f(1-f);
  \ee
  \be
 \label{m37par}
   D_t^7 u= - D_x^{10}u +u(1-u) \LongA - \l^7 f^{(7)}=- f^{(10)} + f(1-f);
  \ee
  \be
 \label{m39par}
   D_t^9 u = D_x^{10}u +u(1-u) \LongA - \l^9 f^{(9)}= f^{(10)} + f(1-f);
 \ee
with the conditions \ef{BC1}.

 \subsection{Results IV: higher-order hyperbolic  PDEs (Section \ref{SHyp1})}

Finally, we consider four purely {\em hyperbolic} higher-order
equations:
 \be
 \label{m32hyp}
   u_{tt}=  D_x^{10}u +u(1-u) \LongA  \l^2 f''= f^{(10)} + f(1-f);
  \ee
\be
 \label{m34hyp}
   u_{tttt}= - D_x^{10}u +u(1-u) \LongA  \l^4 f^{(4)}=- f^{(10)} + f(1-f);
  \ee
  \be
 \label{m36hyp}
   D_t^6 u=  D_x^{10}u +u(1-u) \LongA  \l^6 f^{(6)}= f^{(10)} + f(1-f);
  \ee
  \be
 \label{m38hyp}
   D_t^8 u = - D_x^{10}u +u(1-u) \LongA  \l^8 f^{(8)}= - f^{(10)} + f(1-f);
 \ee
with the conditions \ef{BC1}.
 We also present some TW patterns for the corresponding {\em
 elliptic} PDEs such as
 \be
 \label{m32ell}
   u_{tt}= - D_x^{10}u +u(1-u) \LongA  \l^2 f''= - f^{(10)} + f(1-f);
  \ee
\be
 \label{m34ell}
   u_{tttt}=  D_x^{10}u +u(1-u) \LongA  \l^4 f^{(4)}= f^{(10)} + f(1-f);
  \ee
  \be
 \label{m36ell}
   D_t^6 u= - D_x^{10}u +u(1-u) \LongA  \l^6 f^{(6)}= -f^{(10)} + f(1-f);
  \ee
Clearly, as evolution equations, these lead to unstable (and
ill-posed, in Hadamard's sense, 1906) problems, but as certain
special TW patterns, such solutions make sense, though are highly
oscillatory, as one can expect.

\subsection{Results V:  KPP--(10,11) and KPP--(11,12) (Section \ref{S.11.12})}

These are most exotic KPP--models under consideration. As an
example, we consider two of such PDEs, where the order in $t$
exceeds the order in $x$:
 \be
 \label{10.11}
 D_t^{11} u=- D_x^{10}u + u(1-u) \LongA - \l^{11} f^{(11)}=-
 f^{(10)}+f(1-f),
  \ee
 \be
 \label{11.12}
 D_t^{12} u=- D_x^{11}u + u(1-u) \LongA  \l^{12} f^{(12)}=-
 f^{(11)}+f(1-f),
  \ee
  with conditions \ef{BC1}.

\subsection{Results VI: semilinear eleventh-order in time PDEs (Section \ref{S.1-4.11})}

We end up with the following equations, which were not treated
above and are induced by \ef{10.11}:
\be
 \label{1.11}
 D_t^{11} u=- u_x + u(1-u) \LongA - \l^{11} f^{(11)}=-
 f'+f(1-f),
  \ee
\be
 \label{2.11}
 D_t^{11} u=- u_{xx} + u(1-u) \LongA - \l^{11} f^{(11)}=-
 f''+f(1-f),
  \ee
  \be
 \label{3.11}
 D_t^{11} u=- u_{xxx} + u(1-u) \LongA - \l^{11} f^{(11)}=-
 f'''+f(1-f),
  \ee
  \be
 \label{4.11}
 D_t^{11} u=- u_{xxxx} + u(1-u) \LongA - \l^{11} f^{(11)}=-
 f^{(4)}+f(1-f).
  \ee
Note that, according to {\em a priori} bounds, the models
\ef{2.11} and \ef{4.11} belong to the {\em parabolic} type, while
\ef{1.11} and \ef{4.11} can be treated as {\em dispersion} ones.

\subsection{Results VII: a quasilinear dispersion KPP--(11,1) problem (Section \ref{Sn11})}

As in \cite{GKPPII} for parabolic KPP--problems, we claim that
many present results admit extensions to {\em quasilinear} PDEs.
To this end, we treat the following {\em quasilinear} KPP--(11,1)
(cf. \ef{m31}) problem
 \be
 \label{m31n}
u_t=-D_x^{11}(|u|^n u)+u(1-u)  \LongA - \l f'=-
 (|f|^n f)^{(11)}+f(1-f),
 \ee
 where $n>0$ is a fixed parameter. The main distinguished feature
 of {\em degenerate} equations as in \ef{m31n}
 is that the TW profiles $f(y)$ have finite interface at some
 finite $y_0>0$, so we need first to
 explain
 the local (periodic) structure of such solutions as $y \to
 y_0^-$.

\subsection{The origin of $\log t$-front shift and general goals
of the paper}

In Section \ref{Sdisplog}, using a general third-order in $t$
semilinear model, we show, via a kind of an ``affine centre
subspace expansion", how $\log t$-front shift can appear
 for $t \gg 1$ in convergence to a TW pattern.

Concerning our general classification,
  using the KPP-setting, we refer to such problems
 as to the KPP--$(k,l)$, where $k$ stands for the order of the
 differential operator in $x$ and $l$ for the order of the
 derivative in $t$.

\ssk

 Thus, for several  KPP--$(k,l)$ problems, with $k \ge 3$ and $l \ge 1$, the  main
 questions
  to study, here, as in \cite{GKPPI, GKPPII}, are:

\ssk

 {\bf (I)} {\sc The problem of TW existence:} existence of
 travelling waves by using analytical/numerical methods,

\ssk

{\bf (II)} {\sc The problem of ``maximal" speed:} defining, in a
natural sense, the $\o$-limit set $\o(H)$ of the properly shifter
orbit \ef{TW3}, to discuss whether or not, at least, in some
particular KPP--problems,
 \be
 \label{m42}
  \o(H) = \{f(\cdot;\l): \quad \l \in \Lambda\} \whereA \Lambda
  \subset \re \,\,\, \mbox{is bounded}.
   \ee
 Two main questions arise:
  \be
  \label{m44}
\mbox{what is the maximal speed} \quad \l_{\rm max} = \sup
\Lambda, \quad \mbox{and whether} \quad \Lambda = \{\l_{\rm
max}\}?
 \ee
 Those questions are ``remnants" of the KPP setting. It
 turns out that, unlike in \cite{KPP}, for several parabolic
 \cite{GKPPI, GKPPII} and dispersion models
  \be
  \label{max77}
  \Lambda=(0,\l_{\rm max}), \quad \mbox{so that}
  \quad \l_{\rm max}= \sup \Lambda, \,\,\, \l_{\rm max} \not \in
  \Lambda, \,\,\, \mbox{and}\,\,\,\, \l_{\rm min}=0 \not \in
  \Lambda.
  \ee
However, for some PDEs under consideration, there exist {\em
stationary} patterns corresponding to $\l=0$. Therefore, for such
a variety of semilinear PDEs of different orders, any clear
general conclusions on the above speed sets are illusive, and each
problem might have some individual features.

\ssk

 {\bf (III)} {\sc The $\log t$-shift problem:} studying the stability of the TW $f(y;\l_0)$
 and  derivation of the
 $\log t$-shifting of the moving front in the problem \ef{1.1},
 \ef{1.H},  connected with a kind of an ``affine centre subspace
 behaviour" for the rescaled equation (Section \ref{Sdisplog}).

 Overall, as in \cite{GKPPI, GKPPII}, we expect that the $\log t$-shifting phenomenon is
quite a generic property of many nonlinear KPP-type problems
regardless their particular types.

\ssk

Note that, here, we are not able to solve or even discuss one of
the main problem about the TW velocity $\l_0$ (or velocities? --
the $\o$-limit set of the rescaled orbit might include a connected
continuous curve of profiles $\{f(\cdot,\l), \, \l \in \hat
\Delta\}$), which appear in the PDE setting with the Heaviside
initial data \ef{1.H} (cf. {\bf (II)} above). Numerically, this
would require a full set of hard PDE numerical experiments, which
we do not perform in such a generality. The author believes that
using such a full-scale of both the  ODE and PDE numerics is too
exhaustive and inevitable moves the  research into a pure
numerical area, where some important mathematical aspects of the
KPP--ideology, which comprise the main goal of the present study,
could be lost would be essentially reduced to numerical aspects to
appear for sure\footnote{In other words, the author's position
here as follows: (ODEs) numerics, indeed, helps to understand
various features/aspects of modern nonlinear higher-order PDE
theory, but should be used carefully and balanced, i.e., not
shadowing (and not replacing!) any mathematical
ideas/methods/tools/ideology/etc., which could appear from and
connected with
 treating reliable numerical results (though
 the above is more relevant to \cite{GKPPI} than to the present paper).}.


\section{Dispersion equations of various orders in $t$}
 \label{SDisp1}



As in \cite{GKPPI} for PDEs \ef{E4}, \ef{E6}, and \ef{m1}
($m=4,5$), for new classes of dispersion PDEs \ef{m31}--\ef{m39},
it is convenient to present  numerical results, which directly
show the global structure of such TW profiles to be, at least
partially, justified analytically. A more detailed description of
using the {\tt bv4p4c} solver of the {\tt MatLab}
 is given in \cite[\S~2]{GKPPI}. However,
it is important to note that, as the initial data for further
iterations, we always took the Heaviside function as in \ef{1.H}.
This once more had to help us to converge to a proper ``minimal"
profile (indeed, there are many other TW profiles), though, of
course this was not guaranteed {\em a priori}. We keep this rule
for all other KPP--$(k,l)$ problems of interest, including the
hyperbolic ones \ef{m32hyp}--\ef{m38hyp} and more higher-order
ones, where initial velocity, acceleration, etc. were taken zero.

\begin{figure}
\centering
\includegraphics[scale=0.70]{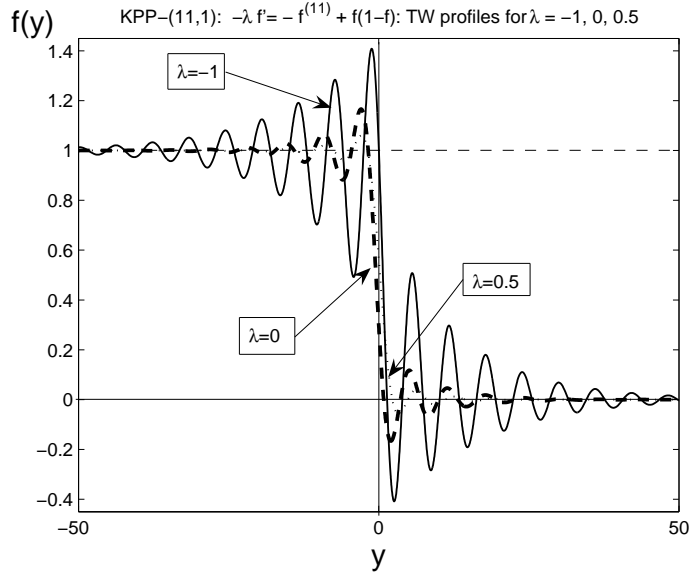}  
\vskip -.3cm
  \caption{A TW profile $f(y)$ satisfying
the ODE in (\ref{m31})  for
 $\l=-1$.}
 \label{F11.1}
\end{figure}

\begin{figure}
\centering
\includegraphics[scale=0.70]{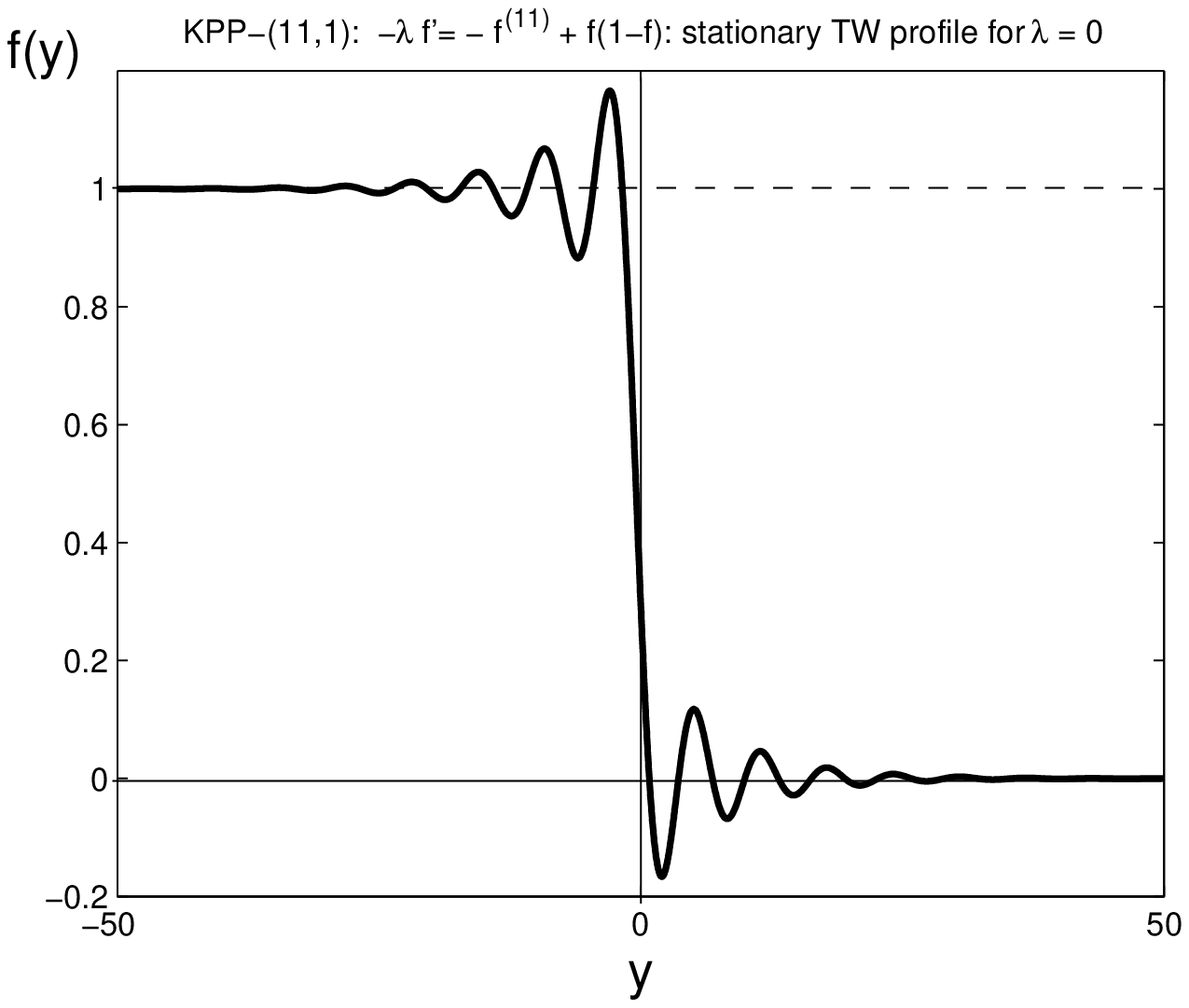}  
\vskip -.3cm
  \caption{The stationary TW profile $f(y)$ satisfying
the ODE in (\ref{m31})  for
 $\l=0$.}
 \label{F11.1St}
\end{figure}

Thus, numerically, we observe existence of TWs for all the PDEs
\ef{m31}--\ef{m39}. In Figure \ref{F11.1}, we show the TW profiles
for the first-order (in $t$) dispersion equation \ef{m31} for the
speeds $\l=-1$, 0.5, and 1. In addition,  a {\em stationary
solution} (denoted by the dash  line in Figure \ref{F11.1}) exists
for $\l=0$, hence, satisfying the ODE (to be referred to later on
a few times)
 \be
 \label{st11.1}
 f^{(11)}=f(1-f) \inB \re, \quad \mbox{plus \,\,(\ref{BC1}).}
 \ee
 Existence (but not uniqueness) of a solution to \ef{st11.1} (and of ones of arbitrary odd order $f^{(2m-1)}$, with any $m \ge 2$)
  is well
 known for a long time; see \cite{Con, Kop}. Uniqueness of such
 stationary solutions of equations like \ef{st11.1} of various odd
 orders was observed numerically; see \cite{GalEng}, where further
 related references were available.
 For convenience and for future use, in Figure \ref{F11.1St}, we
 present the solution of this stationary ($\l=0$) problem \ef{st11.1}.

\ssk

\begin{figure}
\centering
\includegraphics[scale=0.70]{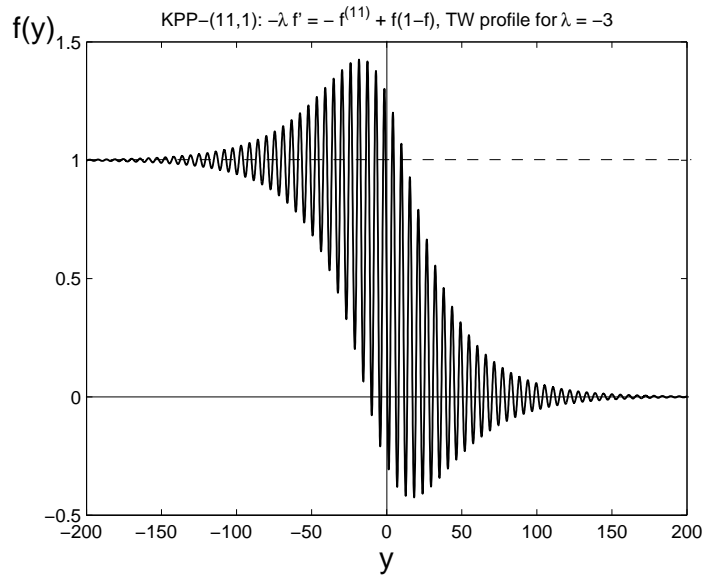}   
\vskip -.3cm
  \caption{An oscillatory TW profile $f(y)$ satisfying
the ODE in (\ref{m31}) and \ef{BC1}  for
 $\l=-3$.}
 \label{F11.1.2}
\end{figure}

Figure \ref{F11.1.2} shows a quite oscillatory TW profile for
$\l=-3$.



\begin{figure}
\centering
\includegraphics[scale=0.70]{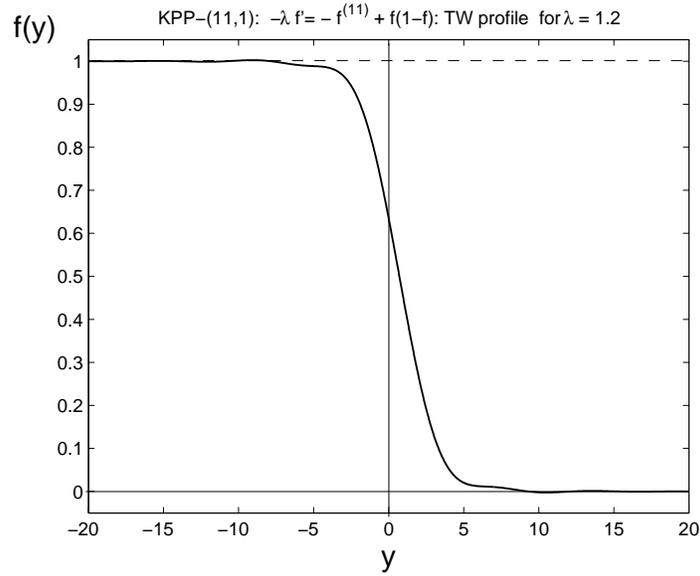}  
\vskip -.3cm
  \caption{A TW profile $f(y)$ satisfying
the ODE in (\ref{m31})  for
 $\l=1.2 \approx \l_{\rm max}^-$.}
 \label{FMax1}
\end{figure}

Further numerical experiments confirm the following $\l$-range of
existence of such TW profiles for \ef{m31} (cf. \ef{SpeedMax}):
 \be
 \label{dispMax}
\mbox{KPP--(11,1):}\quad 1.2   \le \l_{\rm max} <1.3 \, .
    \ee
The TW profile for the existence parameter $\l=1.2$ is shown in
Figure \ref{FMax1}. Numerical integration is performed on a
smaller interval $[-60,60]$, that, in view of small oscillatory
behaviour of $f(y)$ in Figure \ref{FMax1}, cannot cause any
problem (though the length of the interval might affect this
critical speed value). Moreover, it is seen that, as $\l \to
\l_{\rm max}^-$, the profile $f(y)$ tends to be non-oscillatory at
all, so that this $\l_{\rm max}$, characterizes the case when the
characteristic equation of the linearized operator admits complex
roots with vanishing imaginary part. Then this changes the
dimension of the asymptotic bundles
 and makes matching impossible; cf. Theorem 2.1 in \cite[\S~2]{GKPPI}.


\begin{figure}
\centering \subfigure[(\ref{m33})]{
\includegraphics[scale=0.52]{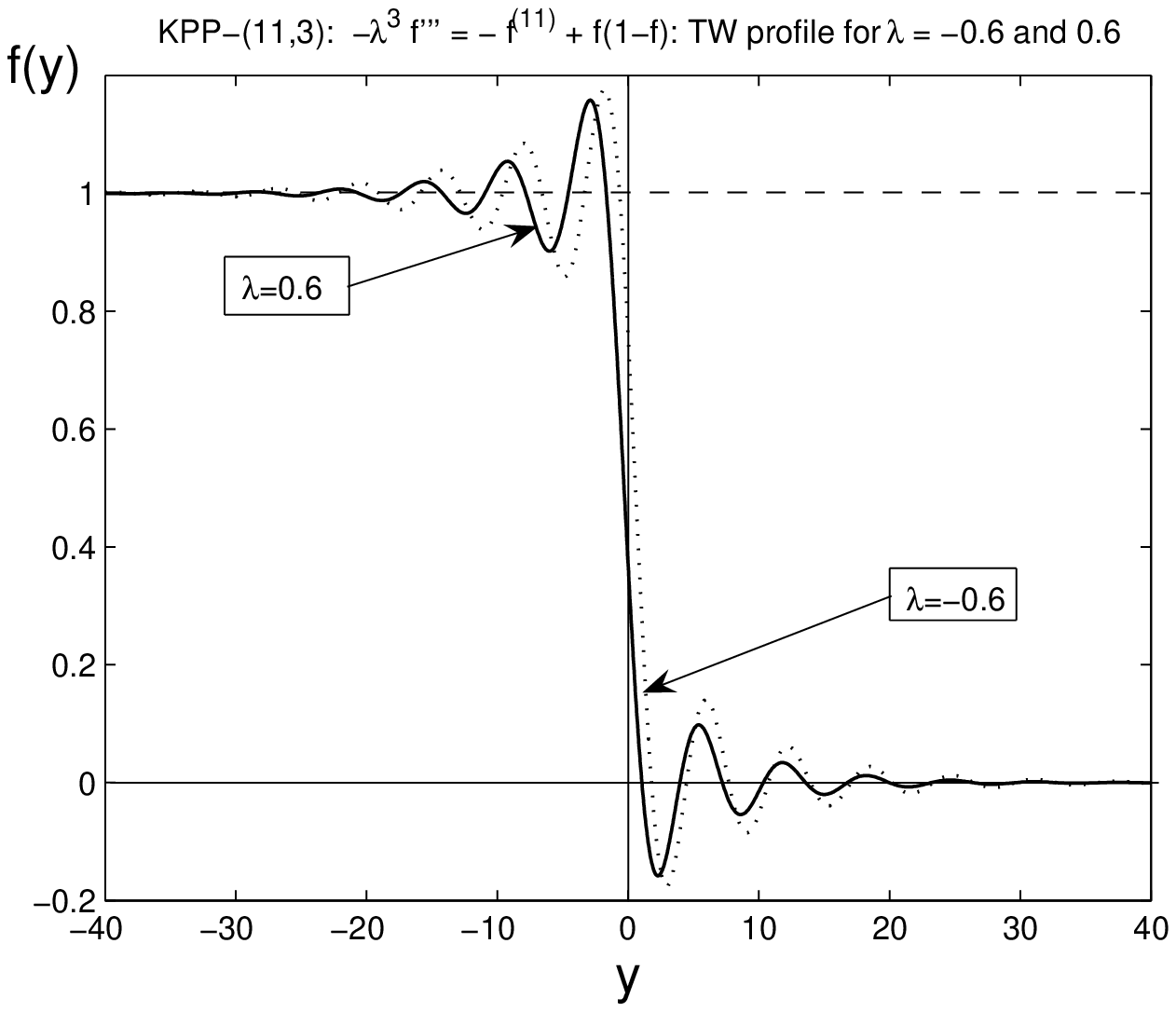}             
} \subfigure[(\ref{m35})]{
\includegraphics[scale=0.52]{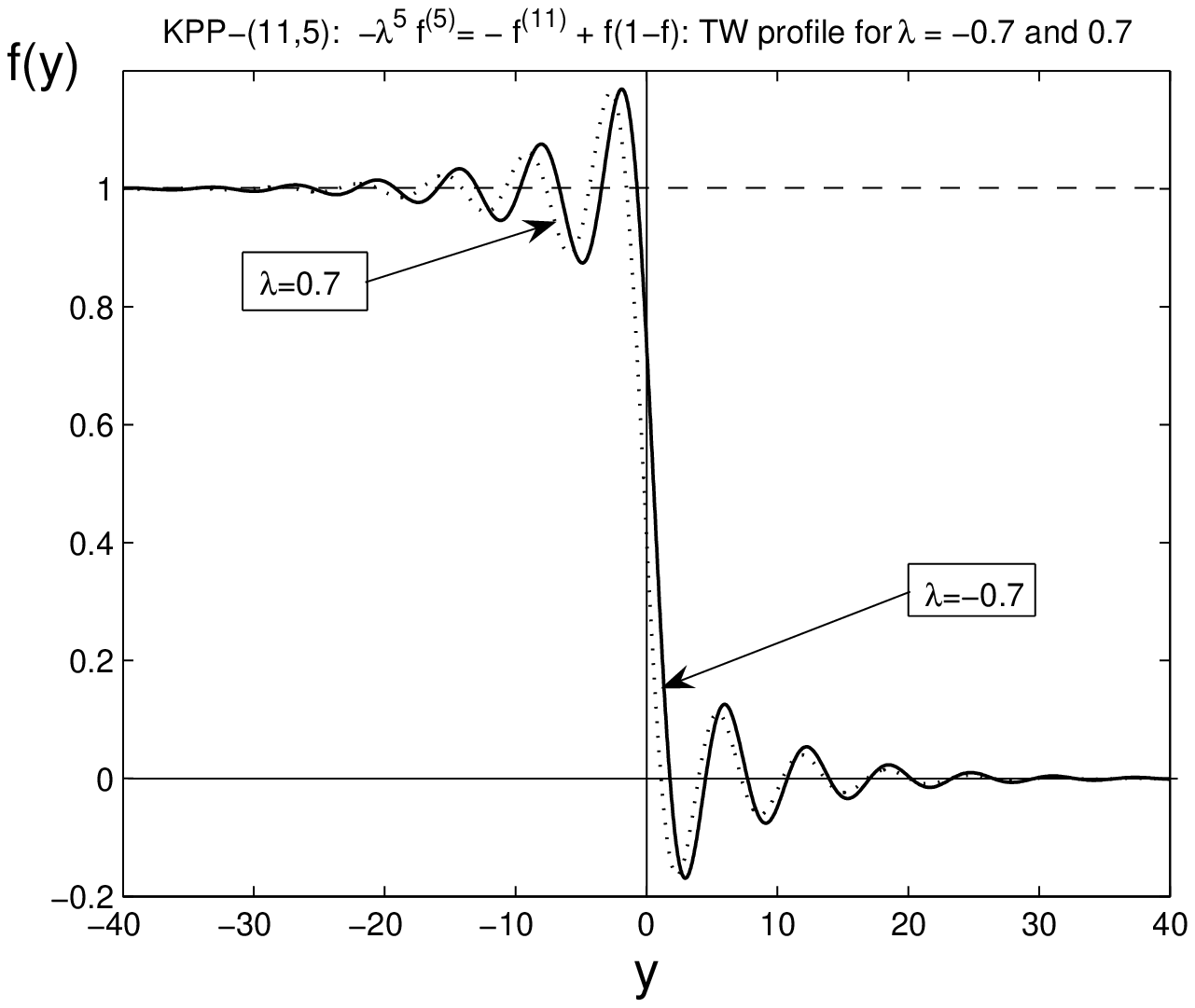}                        
} \subfigure[(\ref{m37})]{
\includegraphics[scale=0.52]{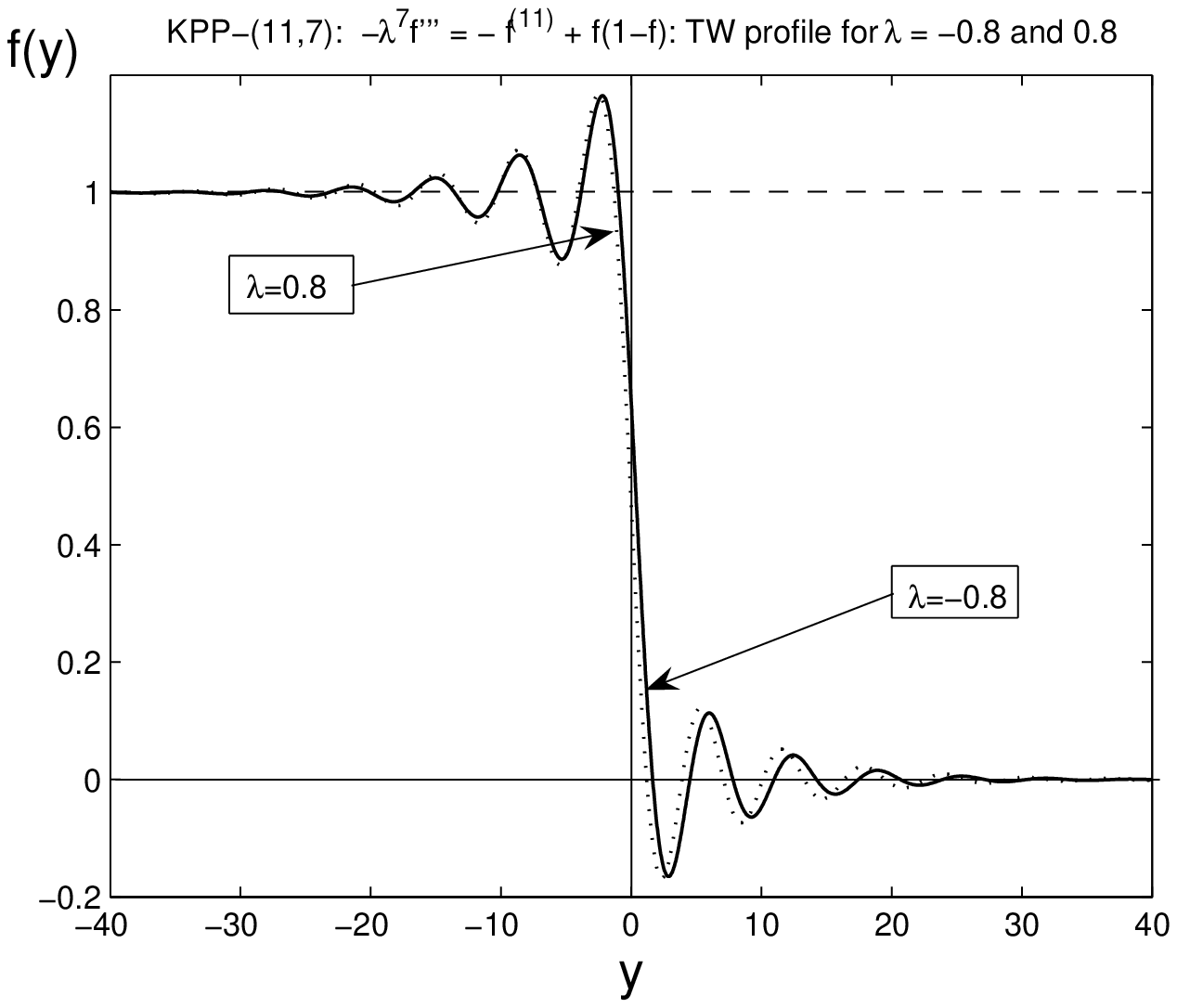}             
} \subfigure[(\ref{m39})]{
\includegraphics[scale=0.52]{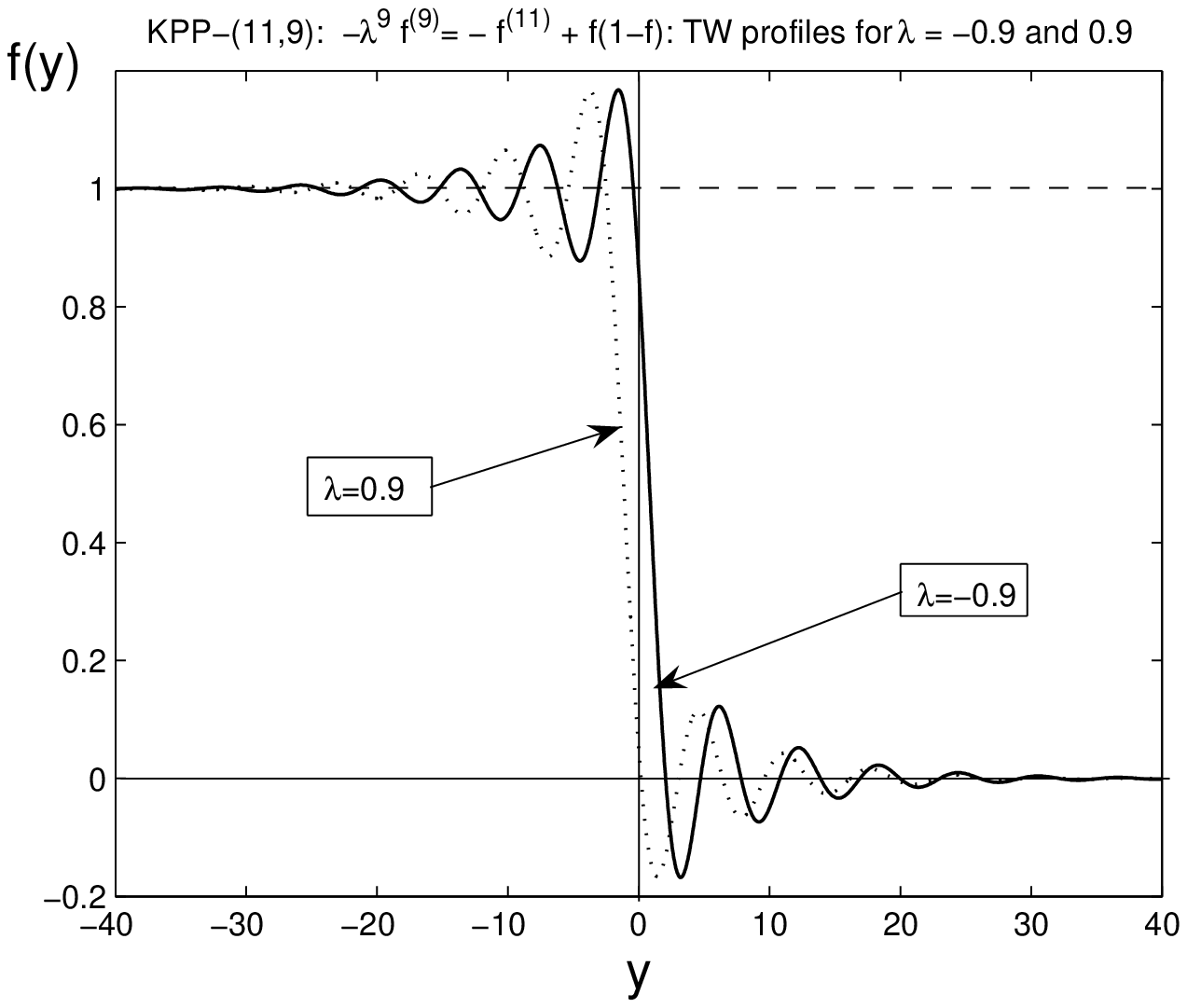}                        
}
 \vskip -.3cm
\caption{\rm\small TW profiles for dispersion equations
\ef{m33}--\ef{m39}.}
 \label{F11.39}
\end{figure}

For equations \ef{m33}--\ef{m39}, TW profiles, with various
negative and positive $\l$'s, are presented in Figure
\ref{F11.39}. In all the cases, there exists also the stationary
profile for $\l=0$ satisfying \ef{st11.1}.

\section{Dispersion-hyperbolic PDEs}
 \label{S.D-H}

\begin{figure}
\centering \subfigure[(\ref{m32dh}), $\l=0,1,0.5$]{
\includegraphics[scale=0.52]{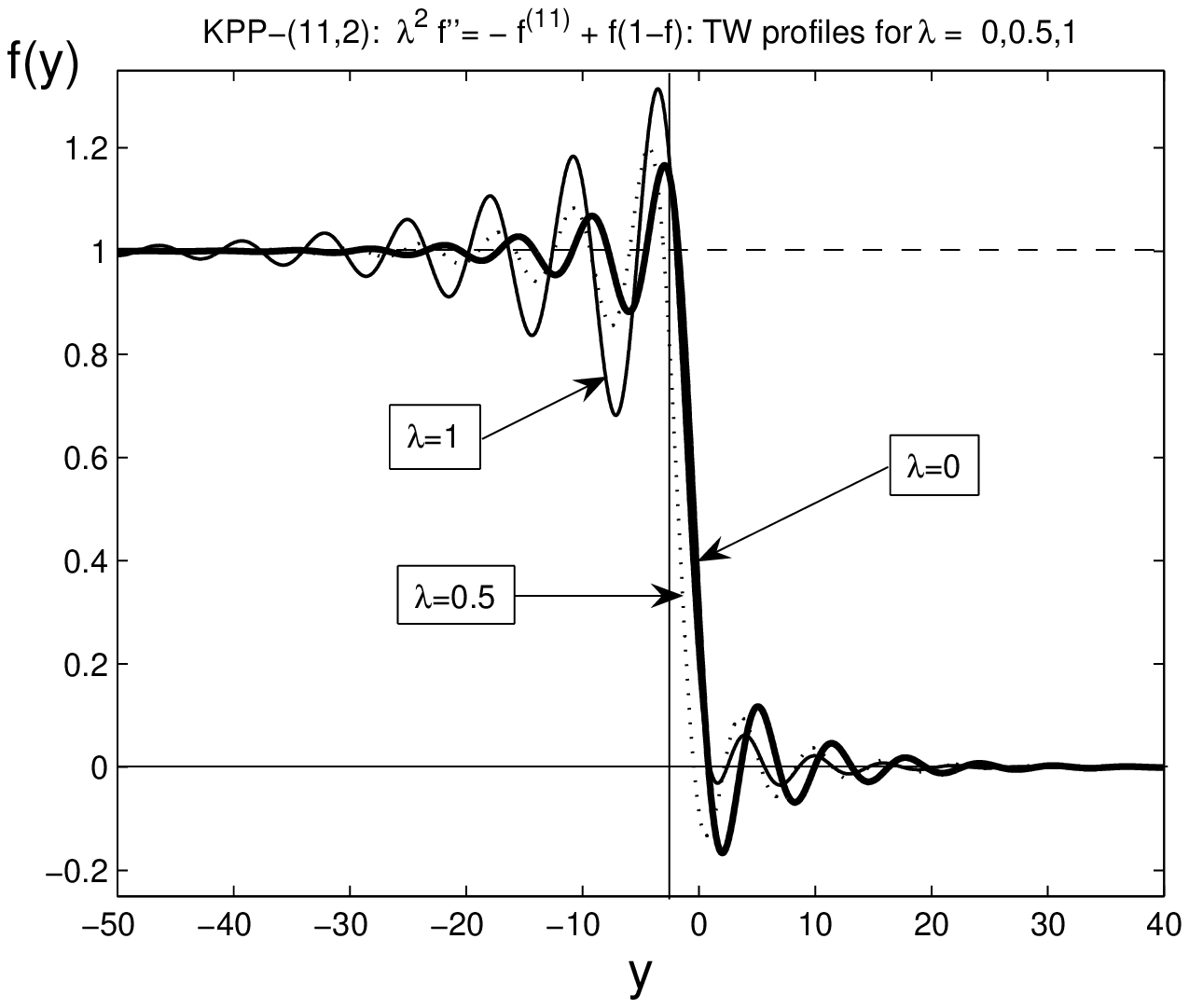}             
} \subfigure[(\ref{m34dh}), $\l=0, 0.6, 1$]{
\includegraphics[scale=0.52]{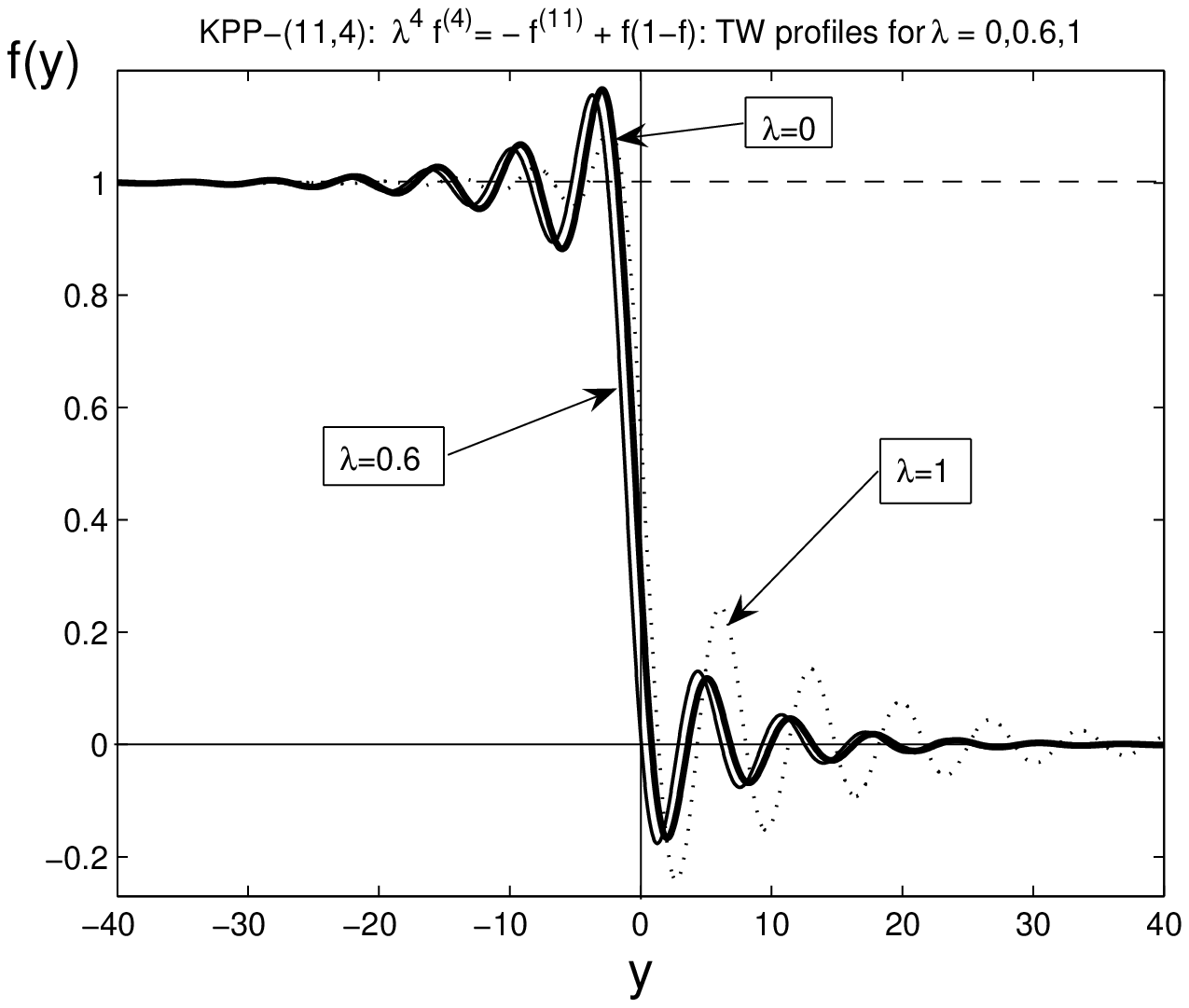}                        
} \subfigure[(\ref{m36dh}), $\l=0.6, 1$]{
\includegraphics[scale=0.52]{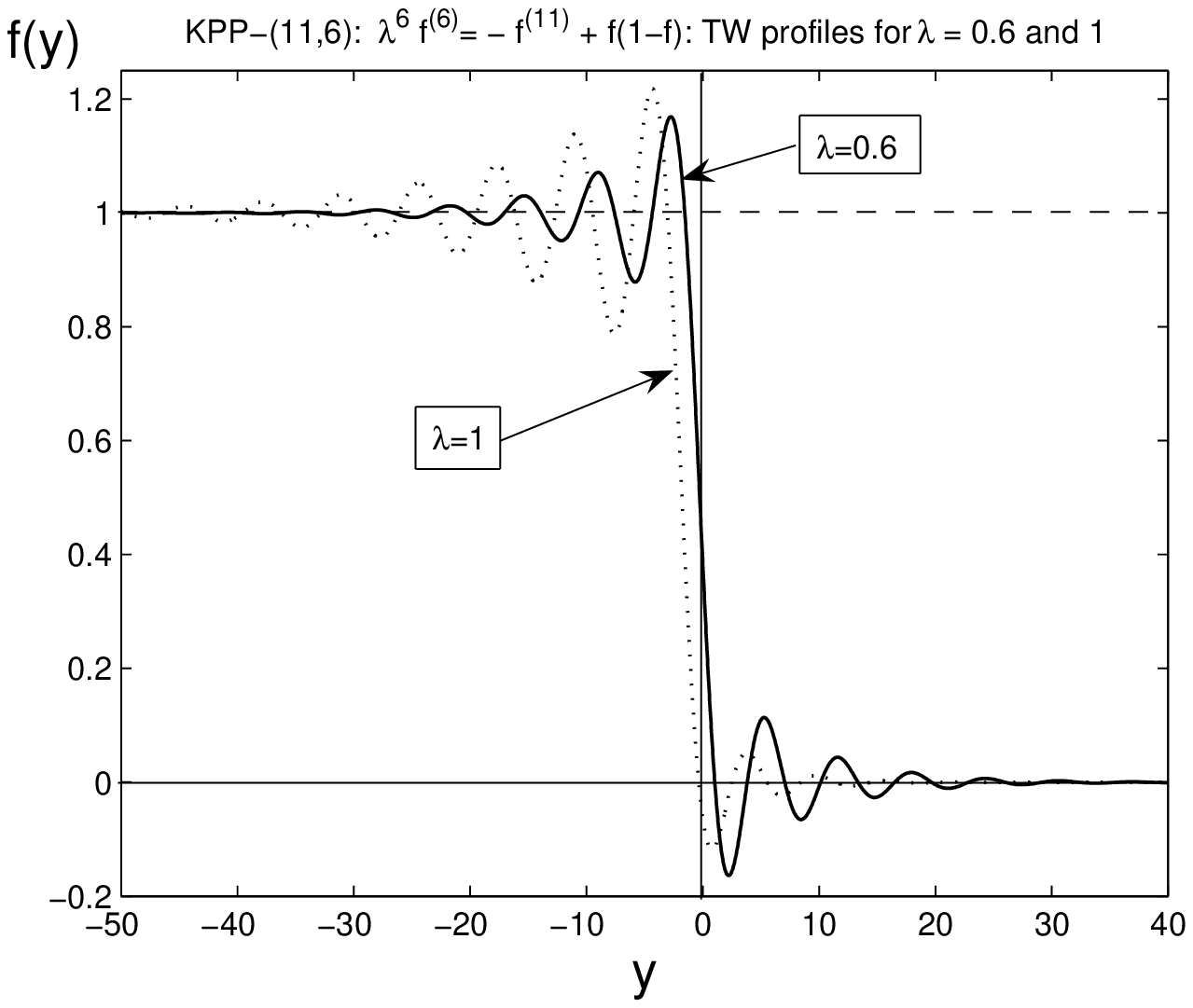}             
} \subfigure[(\ref{m38dh}), $\l=0.06,1,1.04$]{
\includegraphics[scale=0.52]{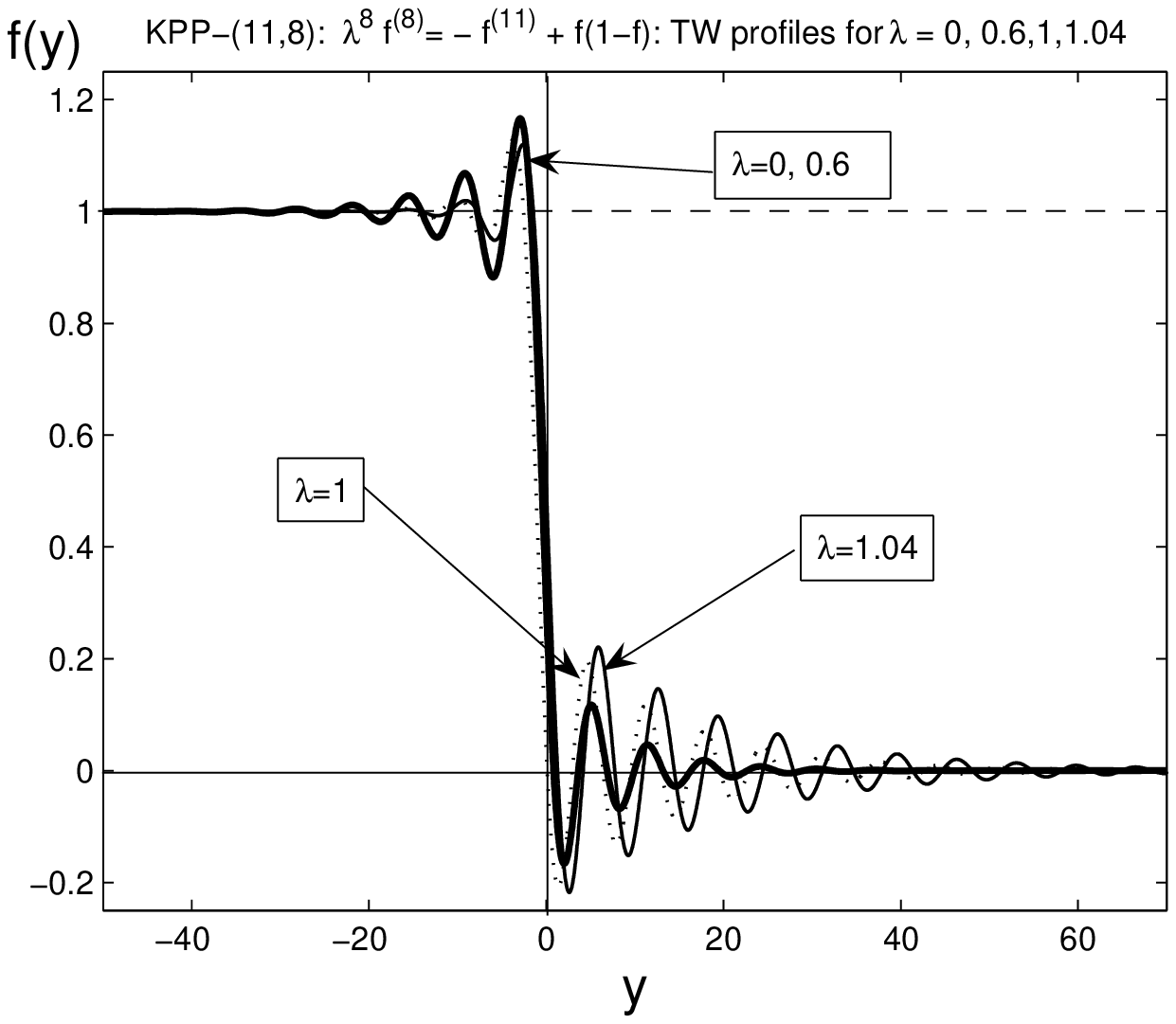}                        
}
 \vskip -.3cm
\caption{\rm\small TW profiles for dispersion-hyperbolic equations
\ef{m32dh}--\ef{m38dh}.}
 \label{F11.dh}
\end{figure}

For equations \ef{m32dh}--\ef{m38dh}, TW profiles, with various
speeds $\l$'s, are presented in Figure \ref{F11.dh}. In the last
case (c), i.e., in the KPP--(11,8), numerics reveal existence of
$\l_{\rm max}$ satisfying
 \be
 \label{ma11}
 1.0443 \le \l_{\rm max} < 1.0445.
  \ee
When approaching $\l_{\rm max}^-$, the TW profiles $f(y)$ remain
essentially oscillatory for $y \gg 1$, so nonexistence for
slightly $\l \ge \l_{\rm max}$ is not related to changing of the
dimension of the  linearized bundle (cf. \cite[\S~2.5]{GKPPI} for
the KPP--(4,1)), but means the impossibility of the corresponding
``nonlinear matching" of such bundles.

\begin{figure}
\centering
\includegraphics[scale=0.70]{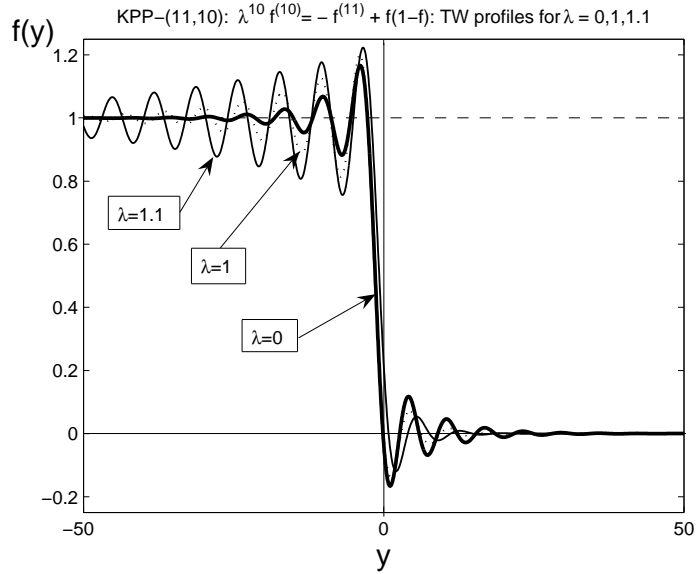}  
\vskip -.3cm
  \caption{TW profiles for the KPP--(11,10) problem \ef{m310dh}.}
 \label{F1110}
\end{figure}

For the most exotic KPP--(11,10) equation \ef{m310dh}, we present
TW profiles for $\l=0$ (the stationary equilibrium), $\l=1$, and
$\l=1.1$ in Figure \ref{F1110}.  Even for $\l$ slightly larger
than 1, the profiles get very oscillatory about the equilibrium 1
 as $y \to -\iy$; cf. $\l=1.1$ in Figure \ref{F1110}.


\section{Further dispersion-parabolic equations}
 \label{S.D-P}

\begin{figure}
\centering \subfigure[(\ref{m33par}), $\l=-0.5$]{
\includegraphics[scale=0.52]{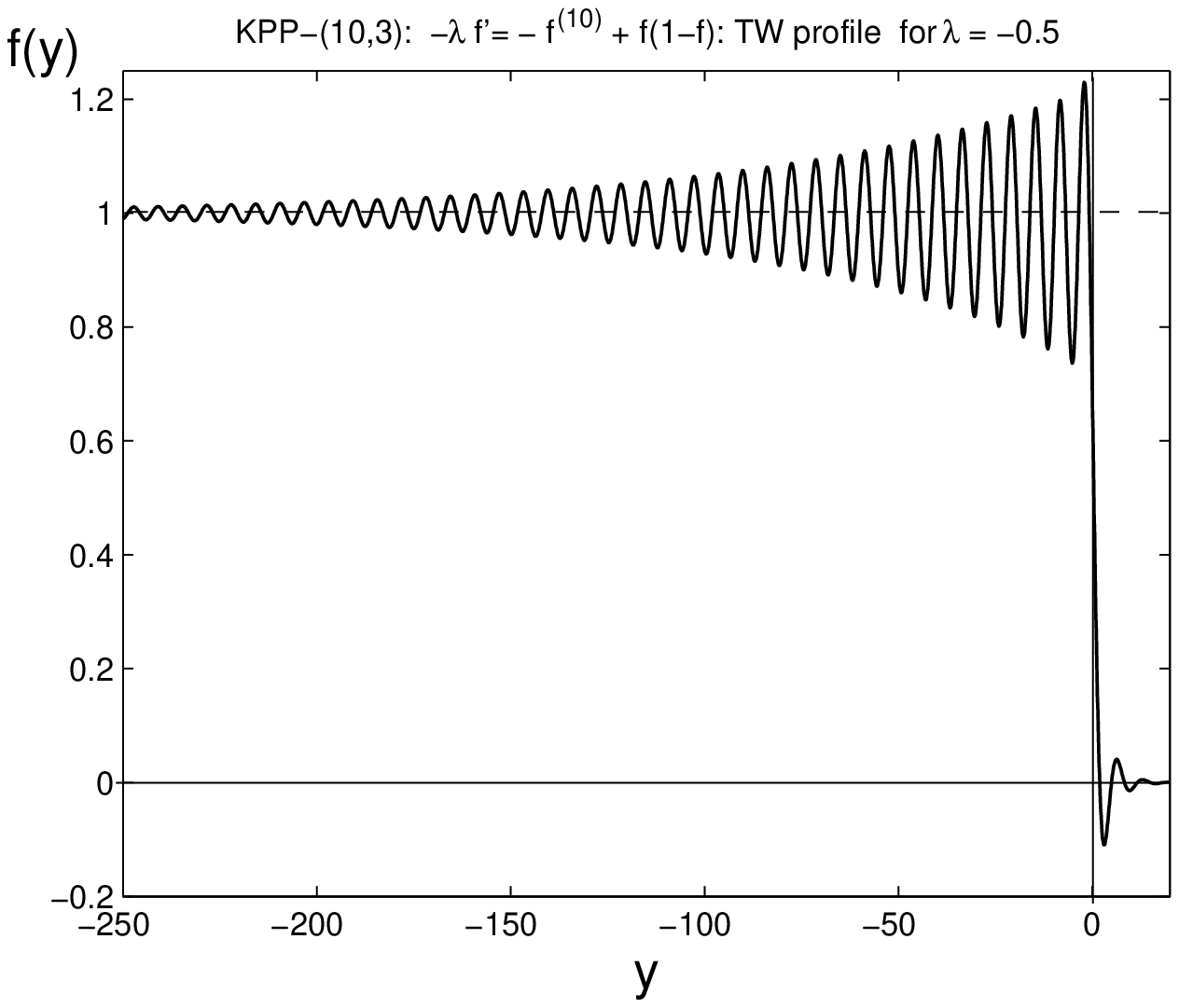}             
} \subfigure[(\ref{m35par}), $\l=1$]{
\includegraphics[scale=0.52]{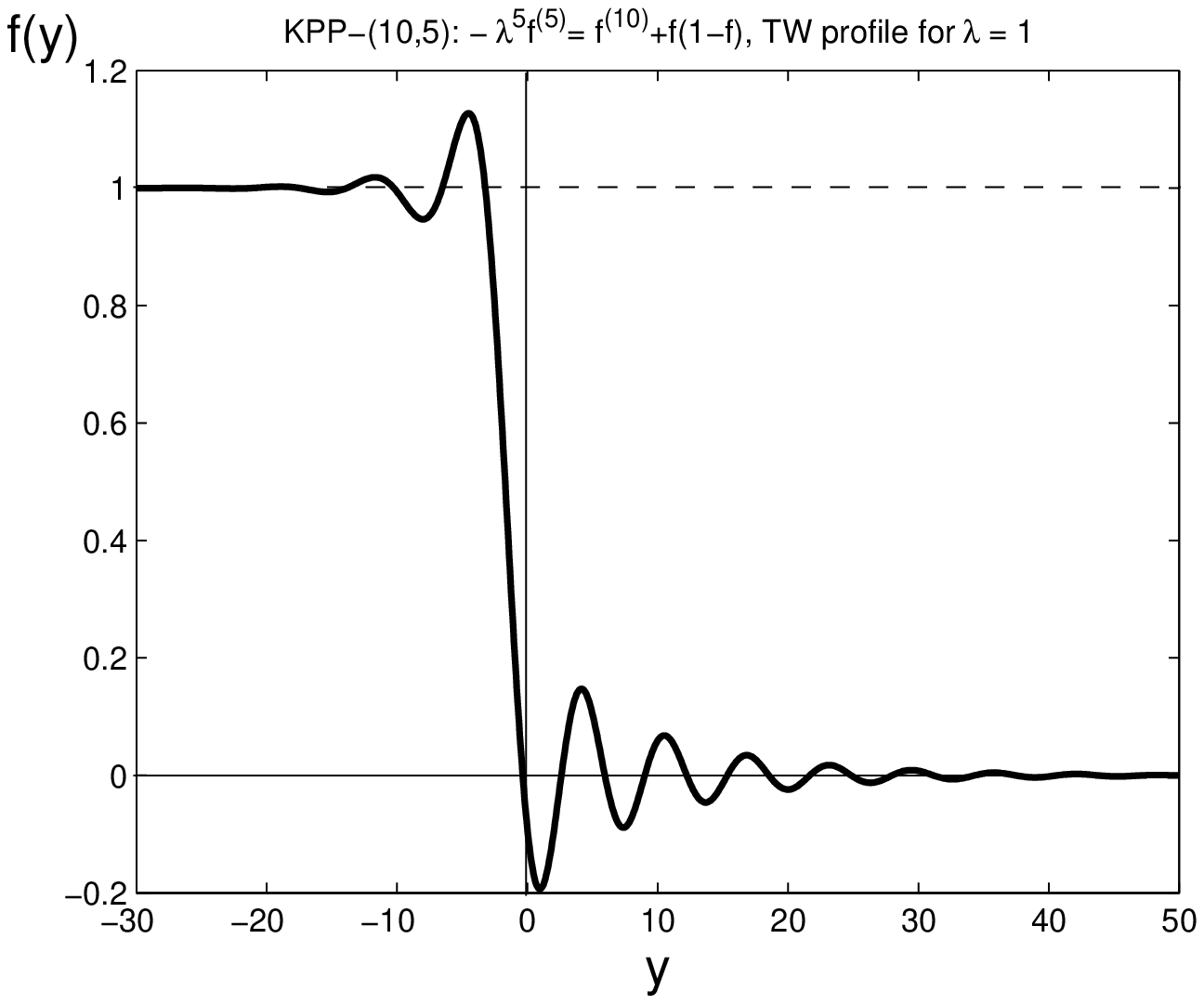}                        
} \subfigure[(\ref{m35par}), $\l=0.5$]{
\includegraphics[scale=0.52]{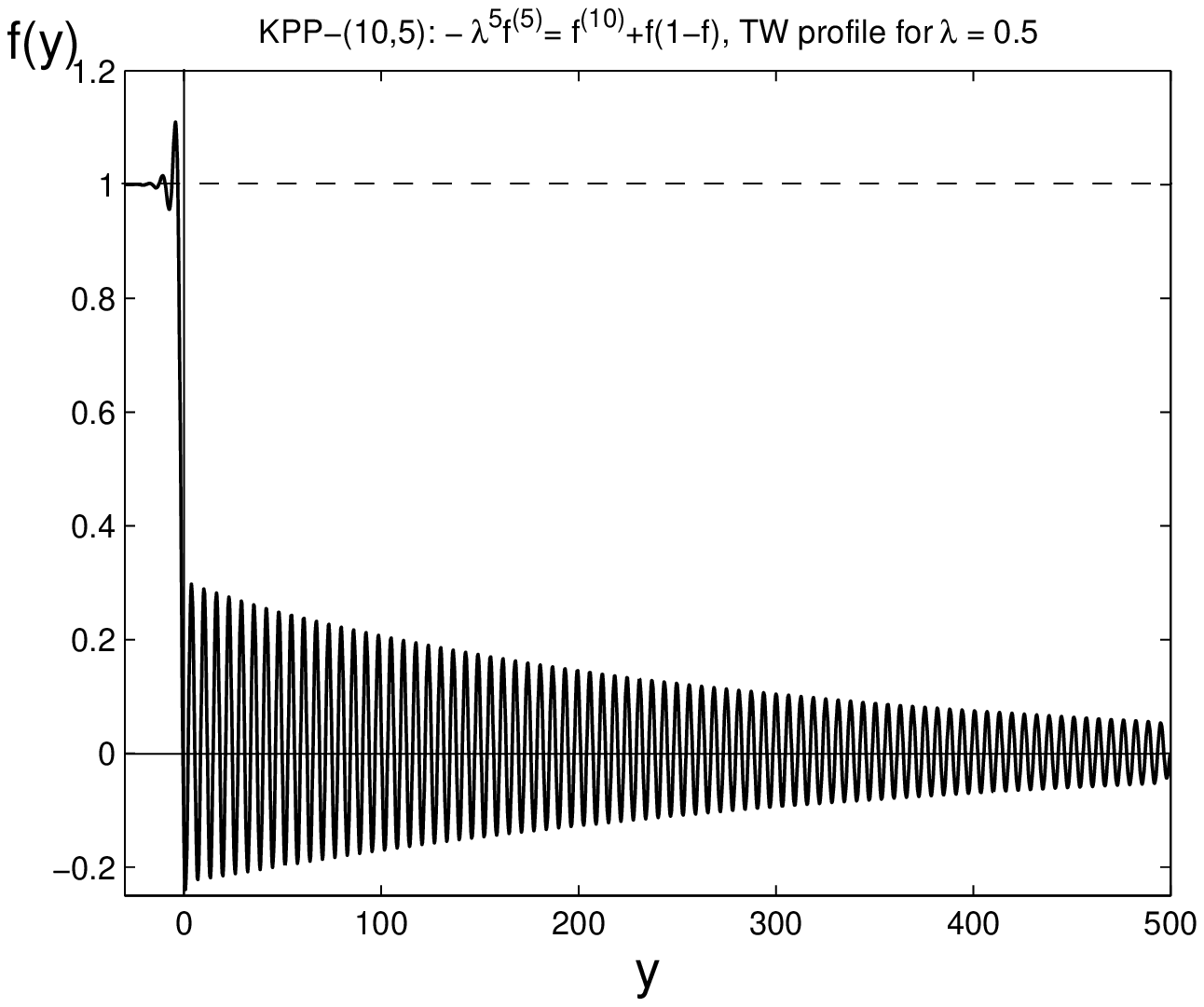}             
} \subfigure[(\ref{m35par}), $\l=1,2,3,4$]{
\includegraphics[scale=0.52]{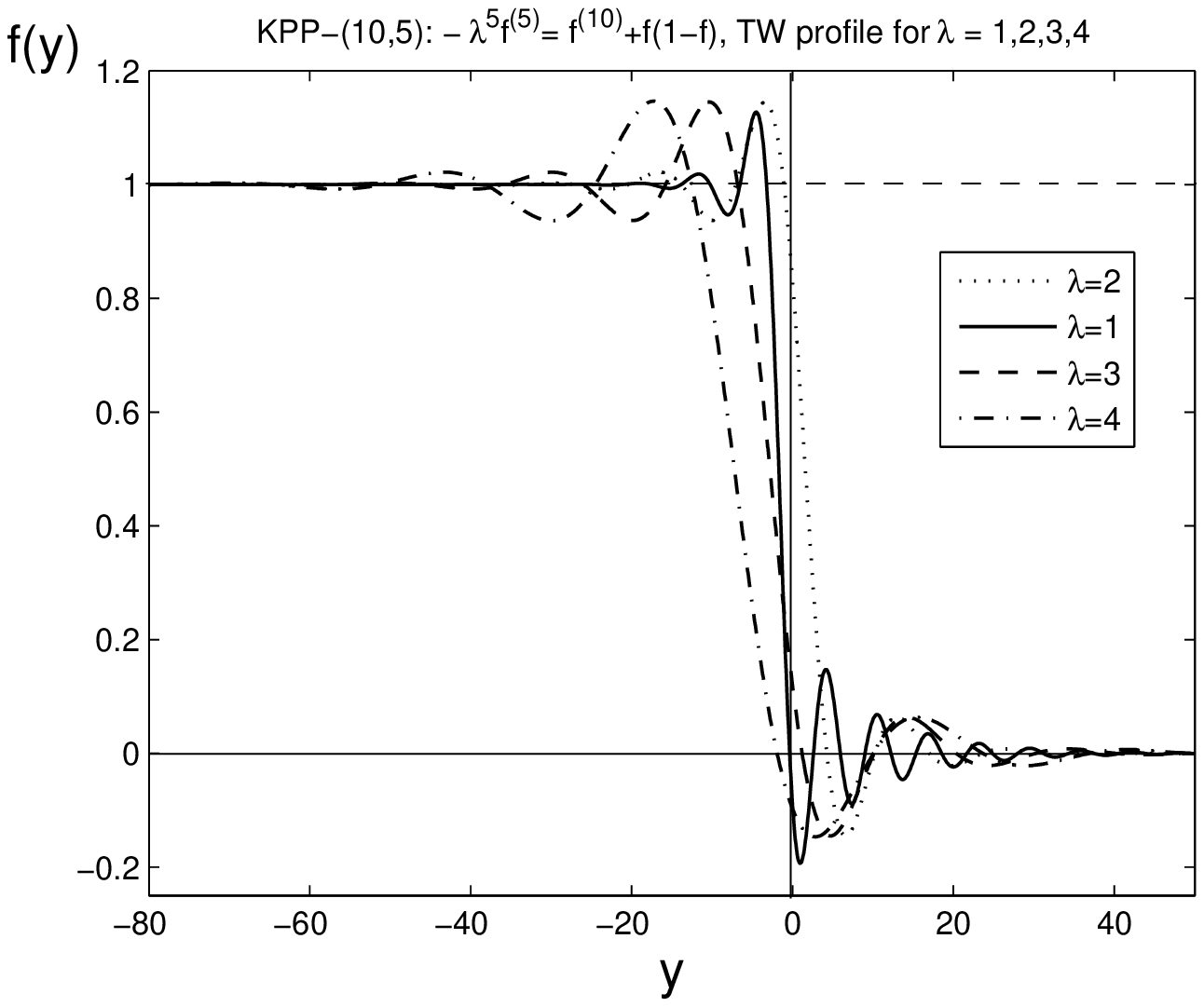}                        
} \subfigure[(\ref{m37par}), $\l=-1$]{
\includegraphics[scale=0.52]{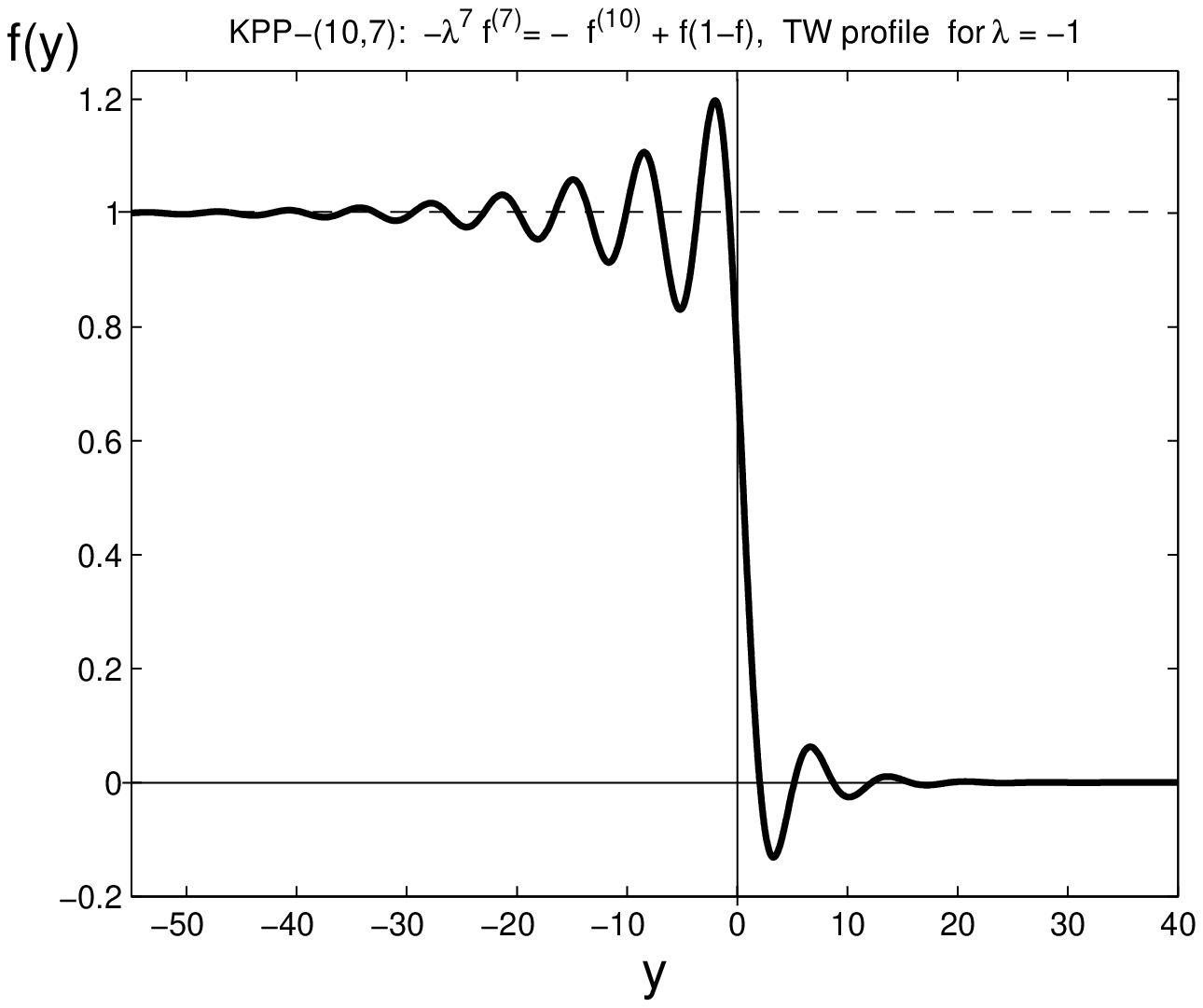}             
} \subfigure[(\ref{m39par}), $\l=1$]{
\includegraphics[scale=0.52]{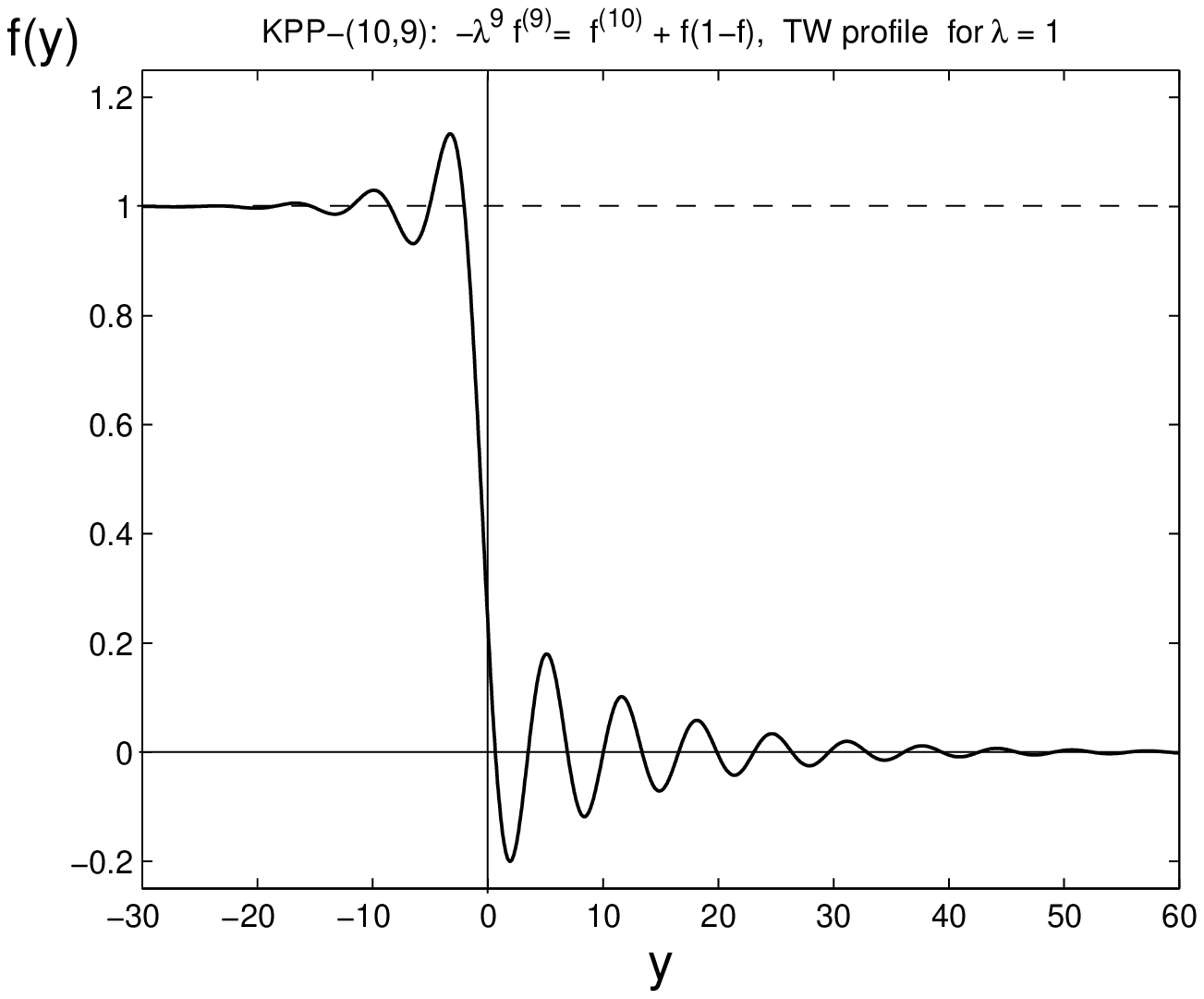}                        
}
 \vskip -.3cm
\caption{\rm\small TW profiles for dispersion-parabolic equations
\ef{m33par}--\ef{m39par}.}
 \label{F10.39}
\end{figure}
\begin{figure}
\centering
\includegraphics[scale=0.70]{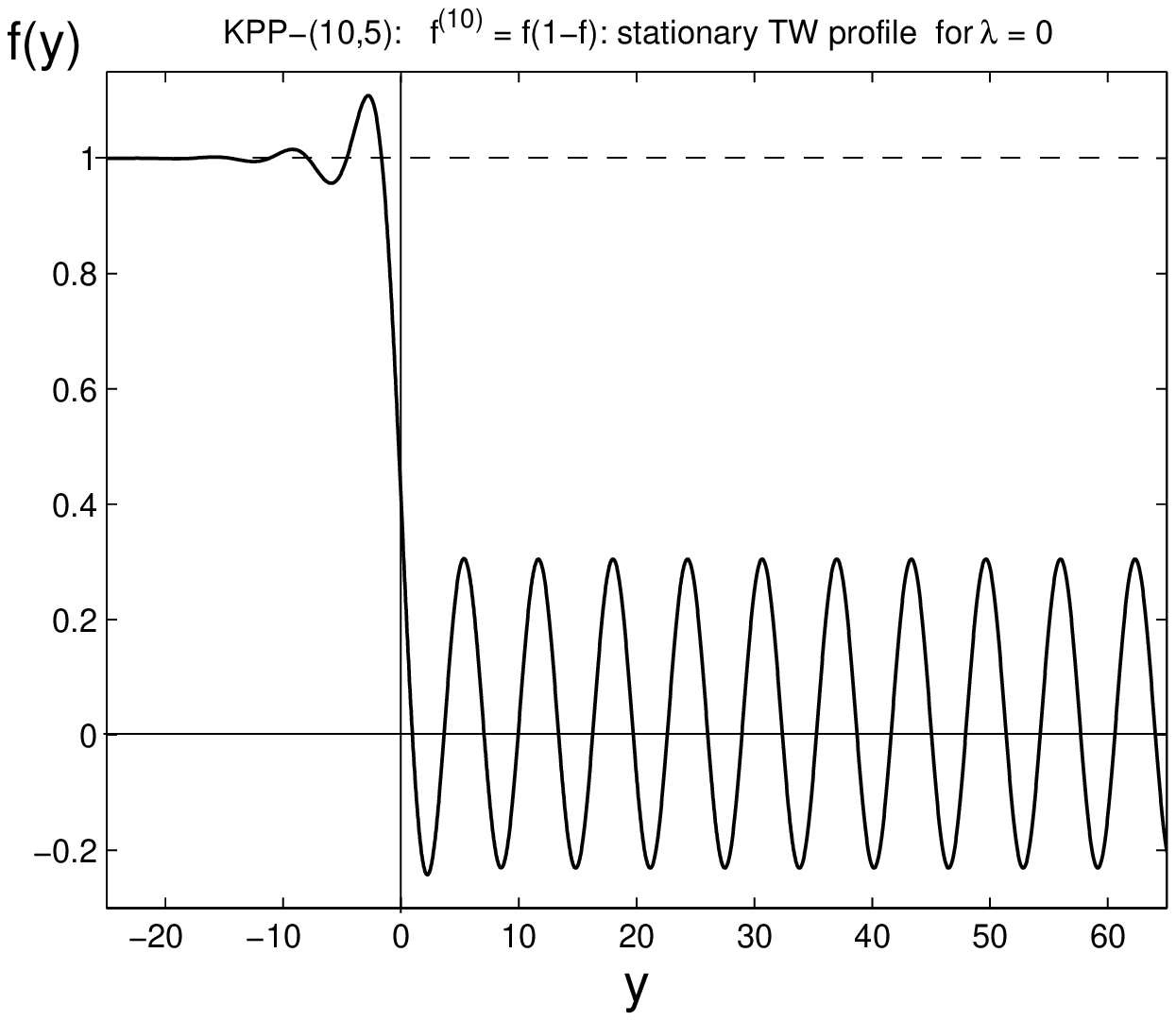}  
\vskip -.3cm
  \caption{A stationary, $\l=0$, solution of \ef{37stat} with a periodic behaviour at $+\iy$.}
 \label{F37Period}
\end{figure}

For equations \ef{m33par}--\ef{m39par}, TW profiles, with various
speeds $\l$'s, are presented in Figure \ref{F10.39}. For equations
\ef{m35par} and \ef{m39par}, for $\l=0$, the stationary equation
 \be
 \label{37stat}
 f^{(10)}=-f(1-f) \inB \re
  \ee
  admits a solution with a periodic behaviour as $y \to +\iy$; see
  Figure \ref{F37Period}.

  Similarly, by symmetry, the stationary equations \ef{m33par} and \ef{m37par}
  admit analogous  stationary profiles that are oscillatory about 1 as $y
  \to -\iy$.


\section{Higher-order hyperbolic equations and  elliptic patterns}
 \label{SHyp1}

\subsection{Hyperbolic equations}

\begin{figure}
\centering \subfigure[(\ref{m32hyp}), $\l=0.5,0.8,1$]{
\includegraphics[scale=0.52]{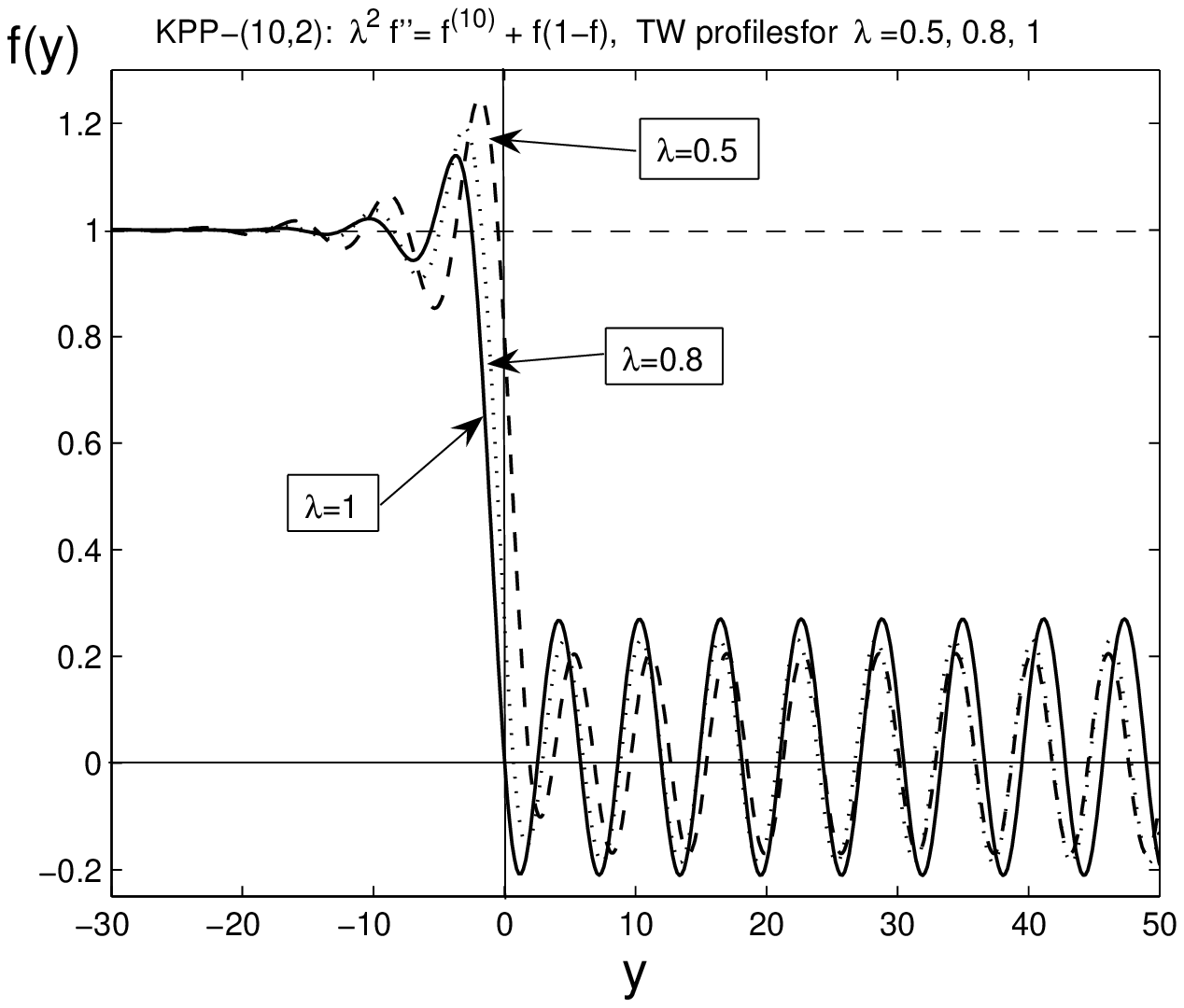}             
} \subfigure[(\ref{m34hyp}), $\l=1$]{
\includegraphics[scale=0.52]{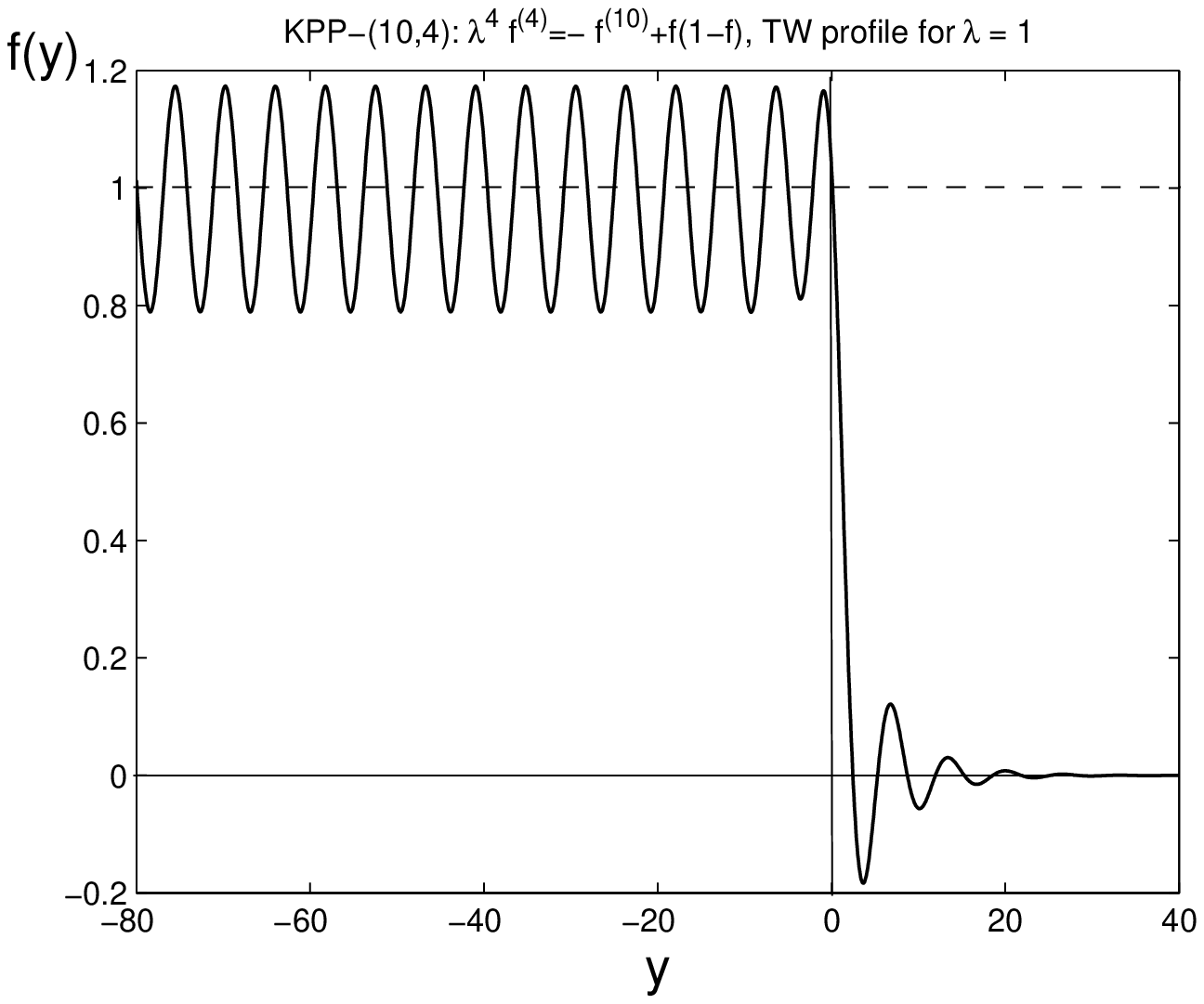}                        
} \subfigure[(\ref{m36hyp}), $\l=-1$]{
\includegraphics[scale=0.52]{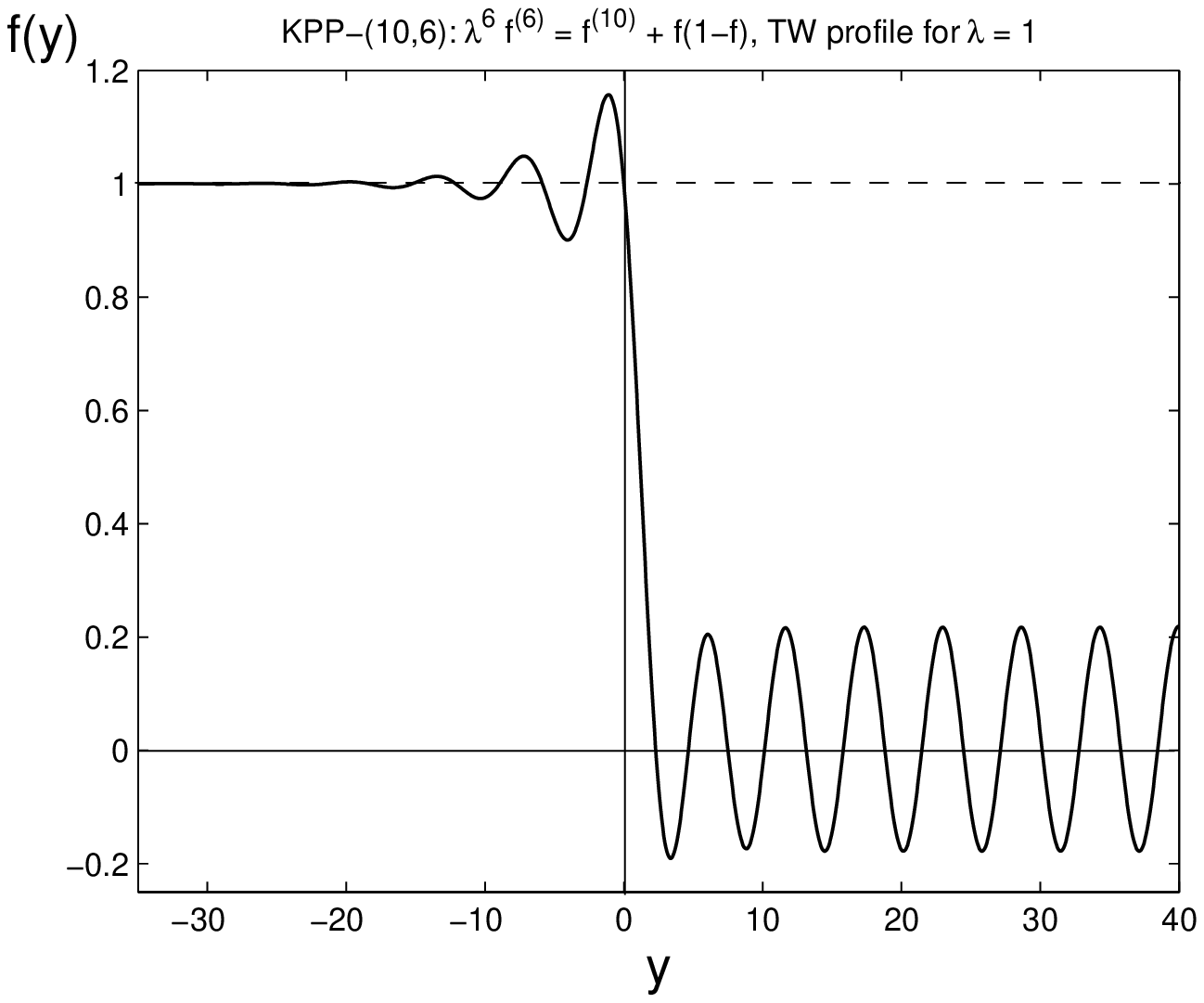}             
} \subfigure[(\ref{m38hyp}), $\l=1,2,3$]{
\includegraphics[scale=0.52]{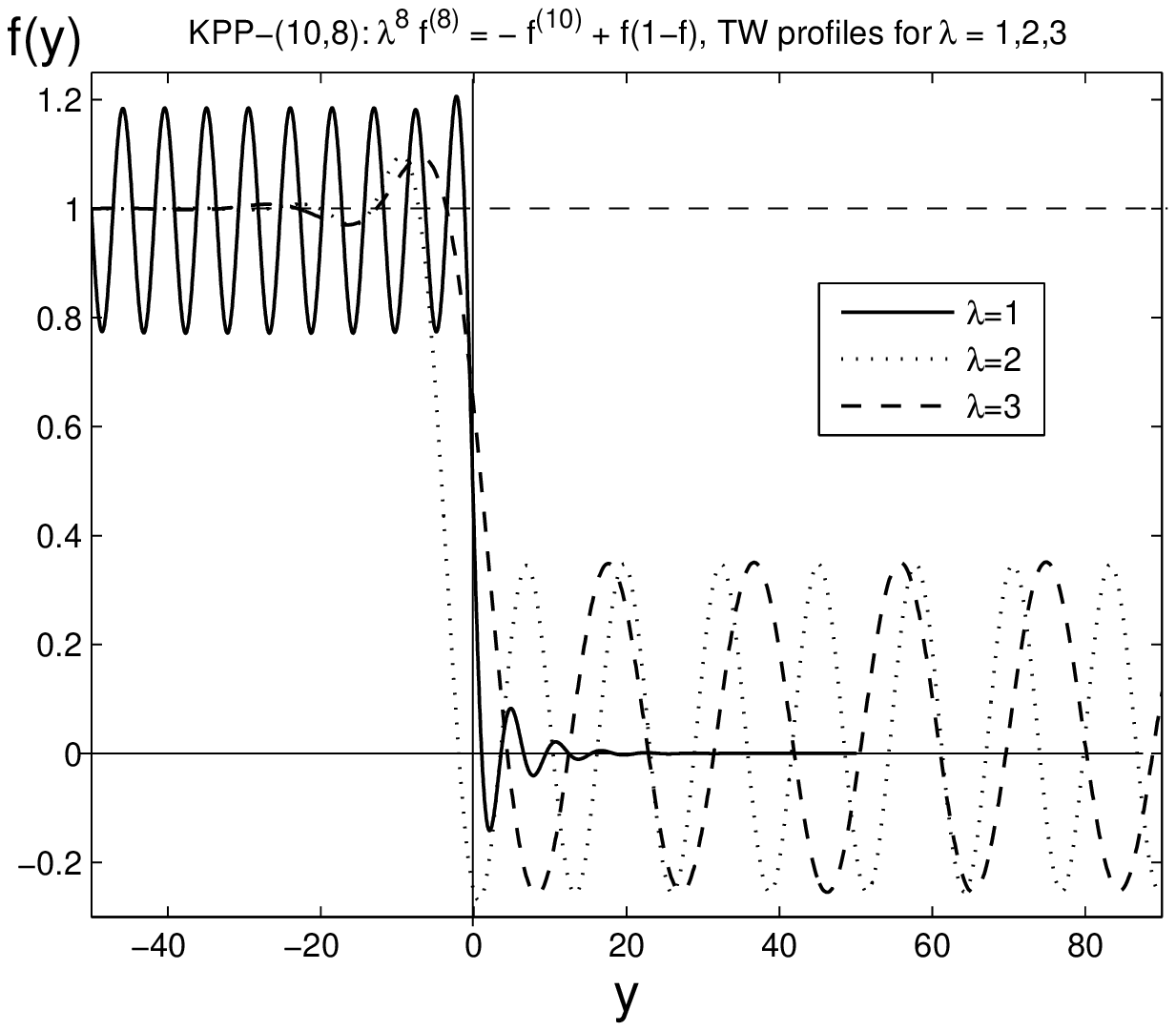}                        
}
 \vskip -.3cm
\caption{\rm\small TW profiles for hyperbolic equations
\ef{m32hyp}--\ef{m38hyp}.}
 \label{F10.hyp}
\end{figure}

For equations \ef{m32hyp}--\ef{m38hyp}, TW profiles, with various
speeds $\l$'s, are presented in Figure \ref{F10.hyp}.

 Note that, for
  such hyperbolic problems, the behaviour  as $y \to +\iy$ often
  becomes
  non-decaying oscillatory (as a feature of hyperbolic flows),
   without a decay or with a very slow decaying
  algebraic envelope. However, as seen from some figures above, several
  TW profiles decay at $+\iy$, and hence satisfy the standard KPP
  setting.

\ssk

  Thus, here,
we  fix operators that are tenth-order in $x$ to avoid any
questions on
 a possibility of a reliable and rigorous local and/or global ODE
 analysis similar to that performed in \cite[\S~2]{GKPPI} for the
 parabolic KPP--4 problem. However, the asymptotic study of these
 ODEs as $y \to +\iy$ can be performed justifying that the dimensions of stable
 bundles as $y \to \pm \iy$ well-correspond to the
  existence
 of TW profiles in the sense of a multi-parameter shooting.
  Nevertheless, any rigorous proof of existence of such $f(y)$
  remains hopeless and represents an open problem, as the existence
  of a {\em heteroclinic path} between the equilibria 0 and
  $\{1,0,...,0\}$ in the tenth-dimensional phase space occurring
  for ODEs in \ef{m32hyp}--\ef{m38hyp}.

Any stability analysis of such TWs in both hyperbolic and
dispersion cases leads to very difficult spectral problems for
pencils of non self-adjoint operators, which we are not going to
treat here; see further ``spectral" references and examples for
hyperbolic/dispersion equations in \cite[\S~8.4]{GalPet2m}.

\subsection{Elliptic patterns}

\begin{figure}
\centering \subfigure[(\ref{m32ell}), $\l=1$]{
\includegraphics[scale=0.52]{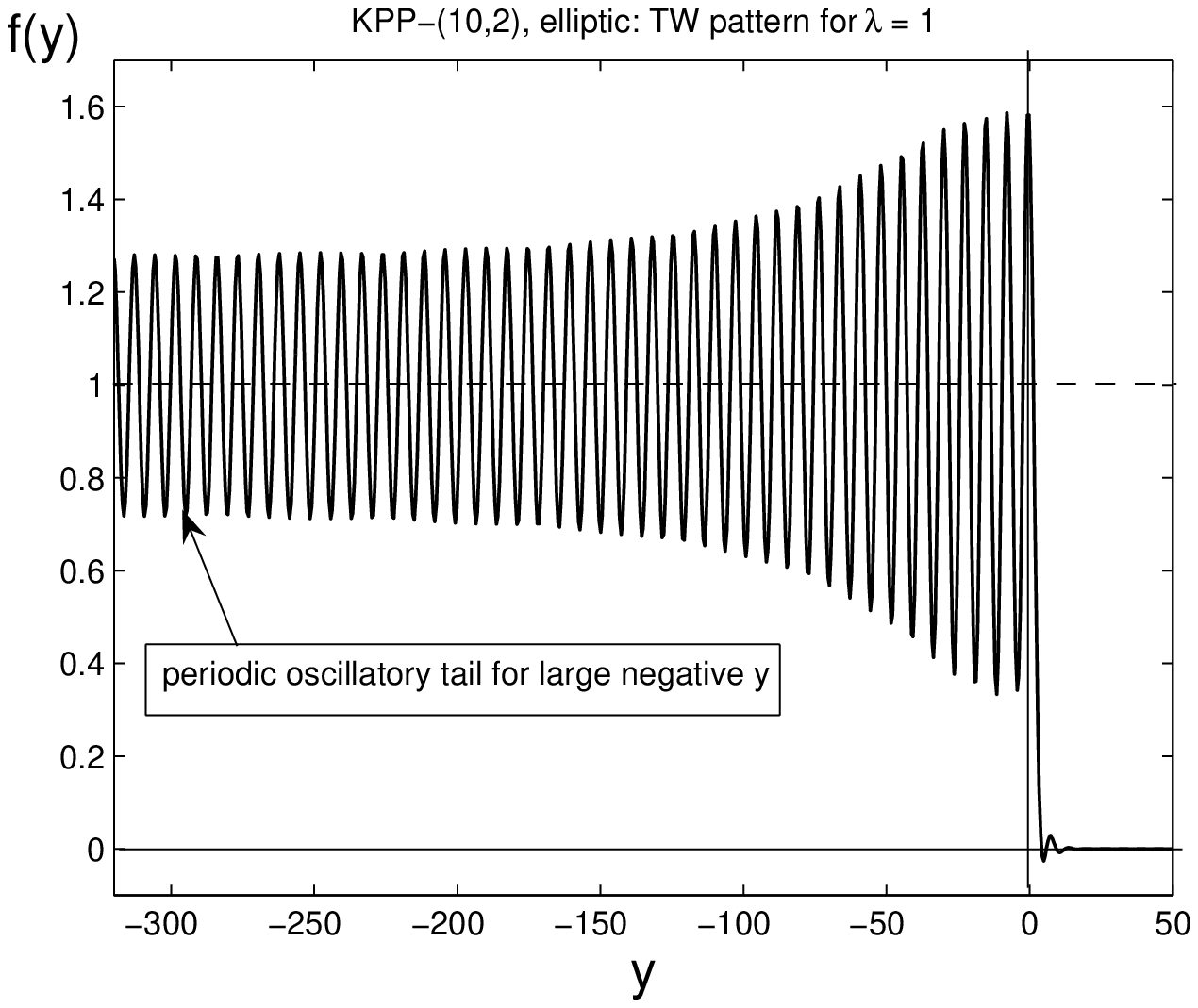}             
} \subfigure[(\ref{m32ell}), $\l=10$]{
\includegraphics[scale=0.52]{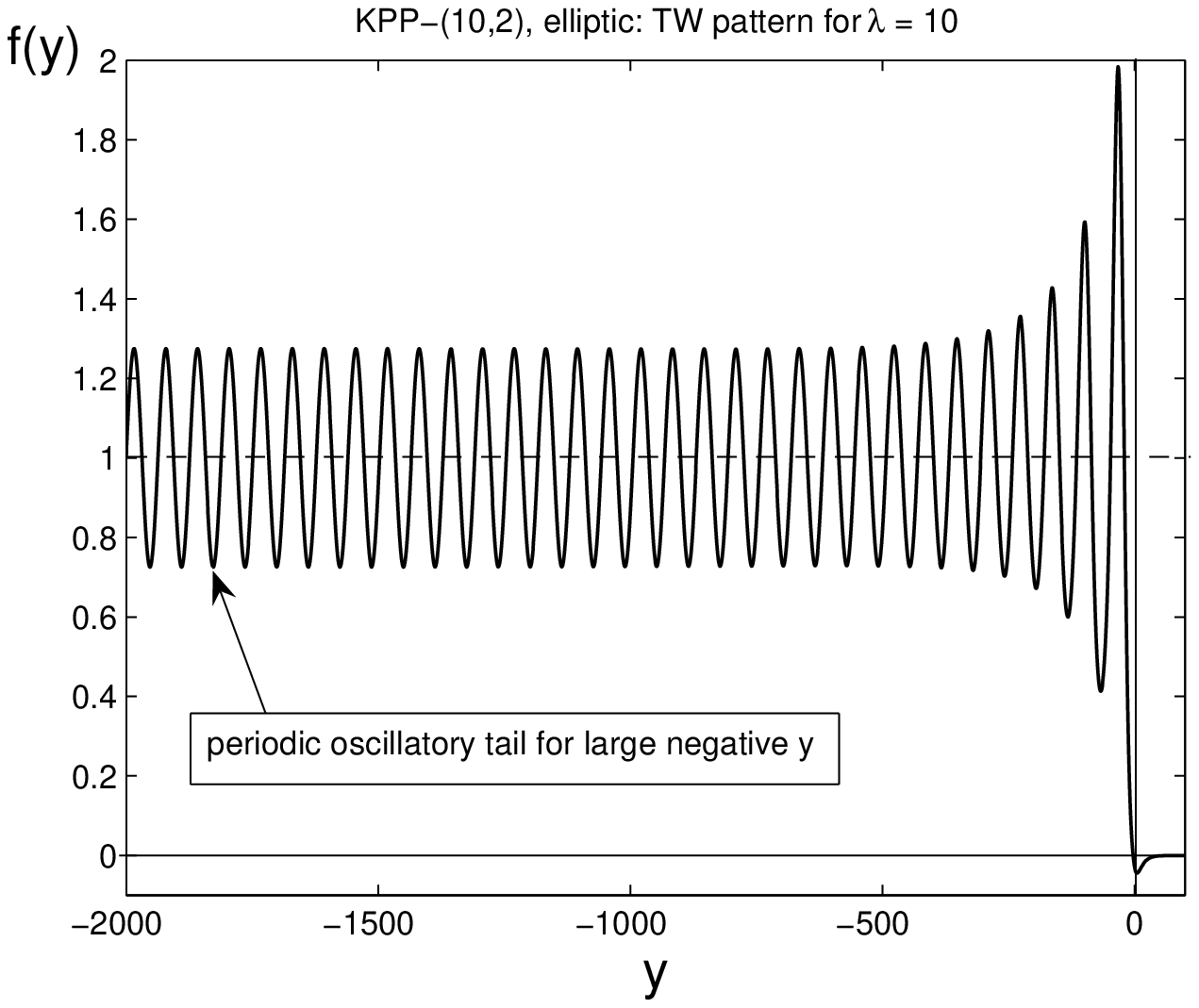}                        
} \subfigure[(\ref{m34ell}), $\l=10$]{
\includegraphics[scale=0.52]{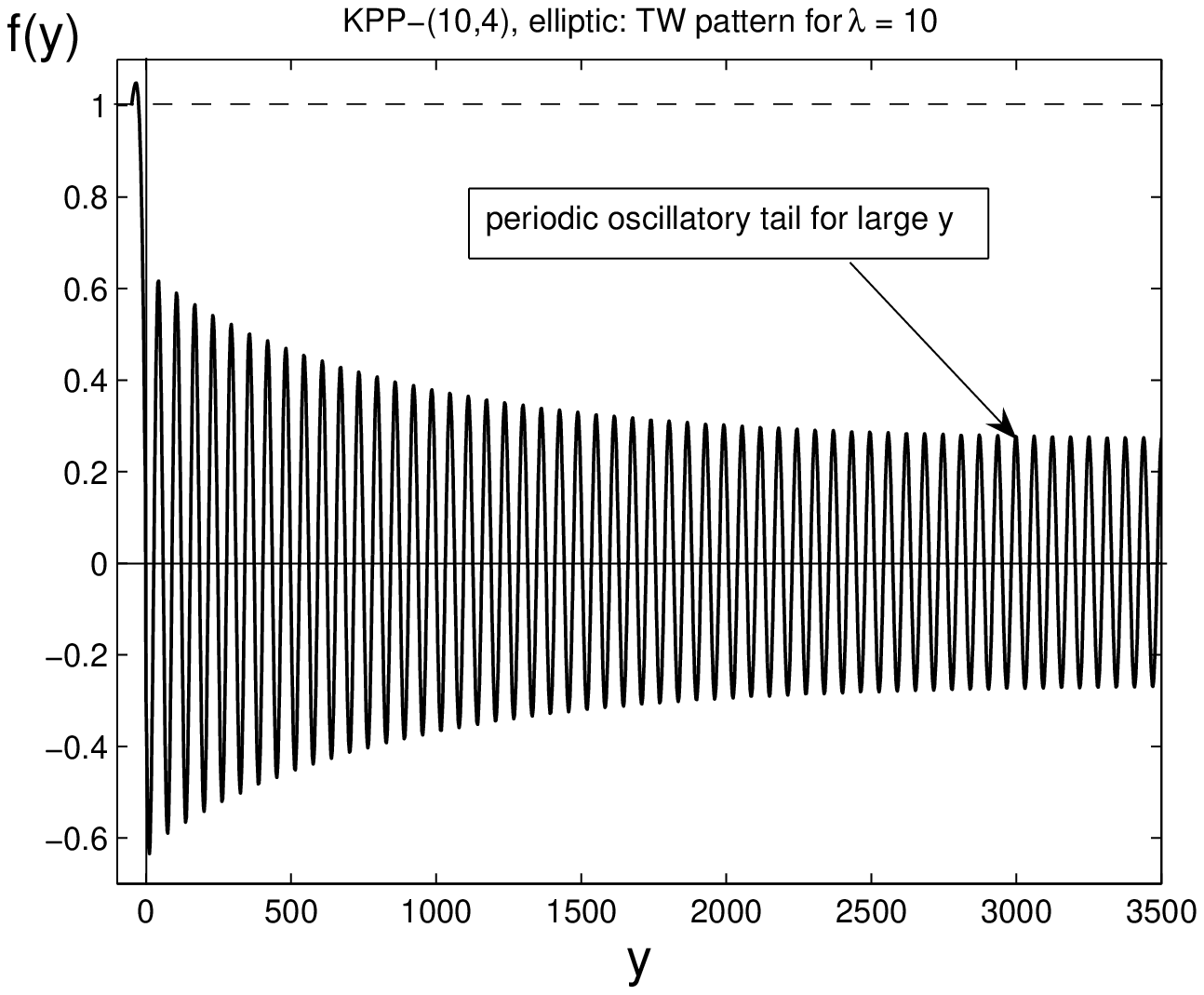}             
} \subfigure[(\ref{m36ell}), $\l=10$]{
\includegraphics[scale=0.52]{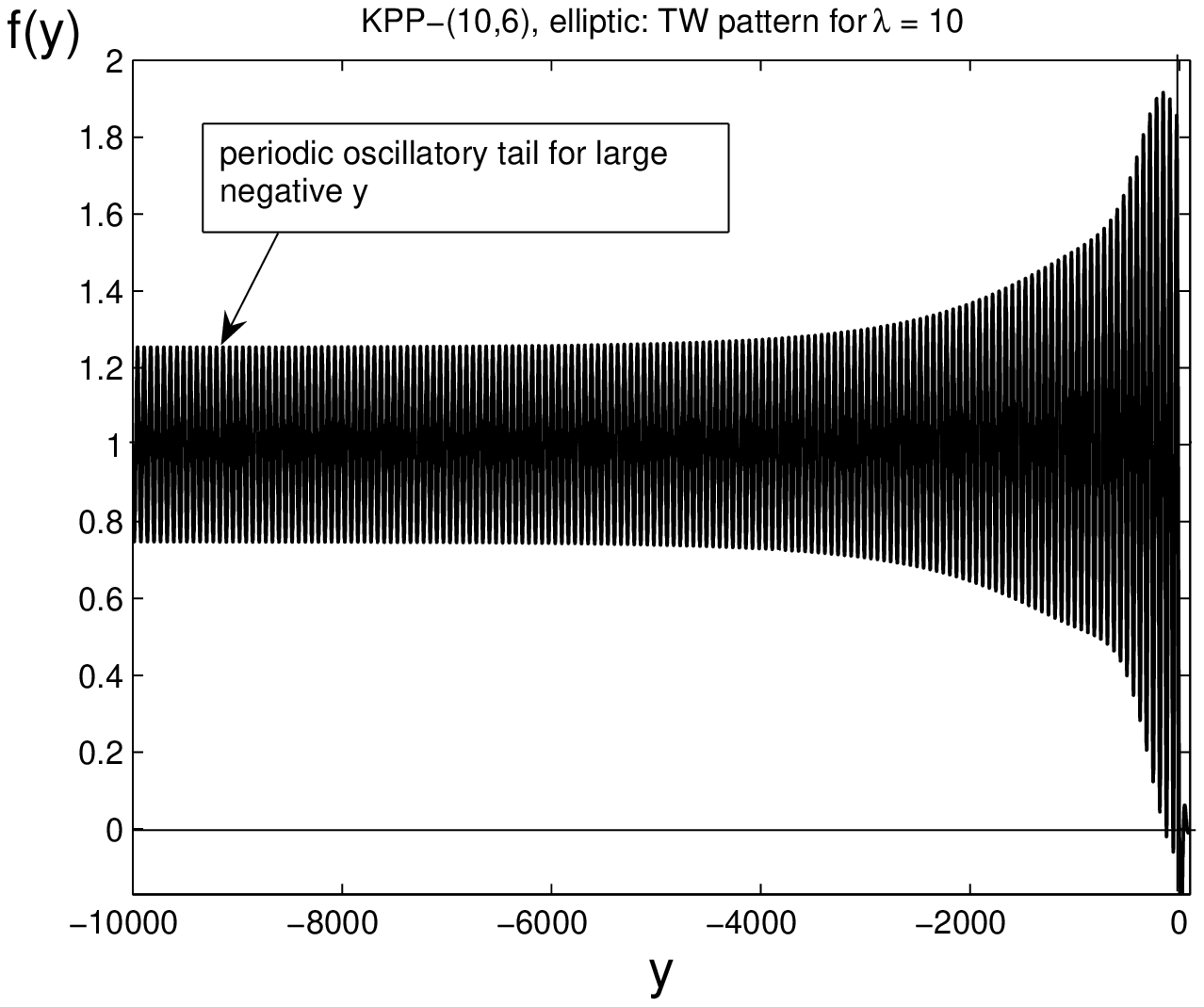}                        
}
 \vskip -.3cm
\caption{\rm\small Oscillatory one-sided ``periodic" TW patterns
for elliptic equations \ef{m32ell}--\ef{m36ell}.}
 \label{F1.Ell}
\end{figure}

Though the corresponding elliptic equations
\ef{m32ell}--\ef{m36ell} do not admit any evolution setting, in
Figure \ref{F1.Ell}, we present elliptic patterns, which are
always highly oscillatory and ``almost periodic" either as $y \to
-\iy$ or $y \to + \iy$.

\section{KPP--(10,11) and KPP--(11,12)}
 \label{S.11.12}

\begin{figure}
\centering \subfigure[(\ref{10.11}), $\l=0$]{
\includegraphics[scale=0.52]{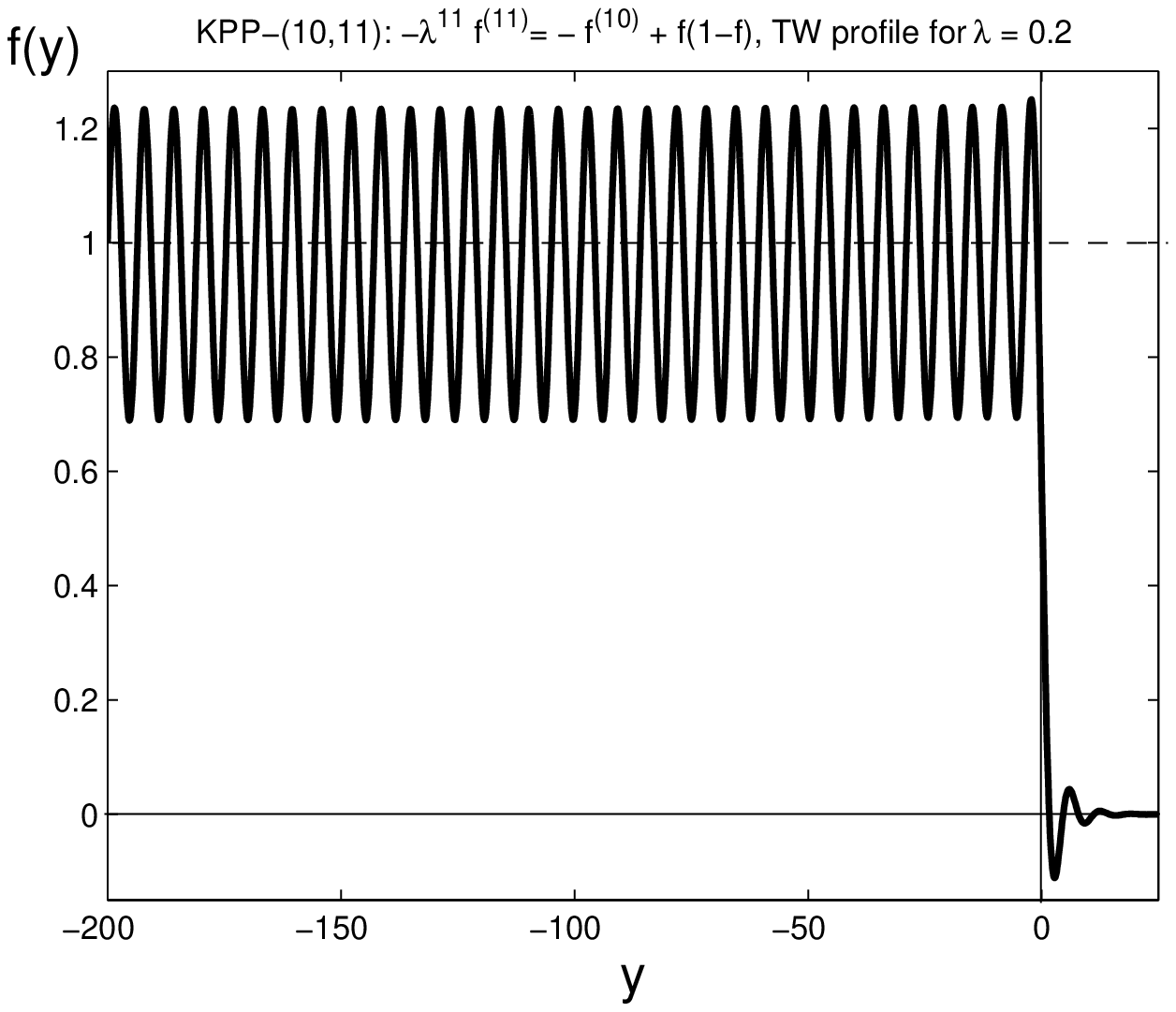}             
} \subfigure[(\ref{10.11}), $\l=-1,-1.1$]{
\includegraphics[scale=0.52]{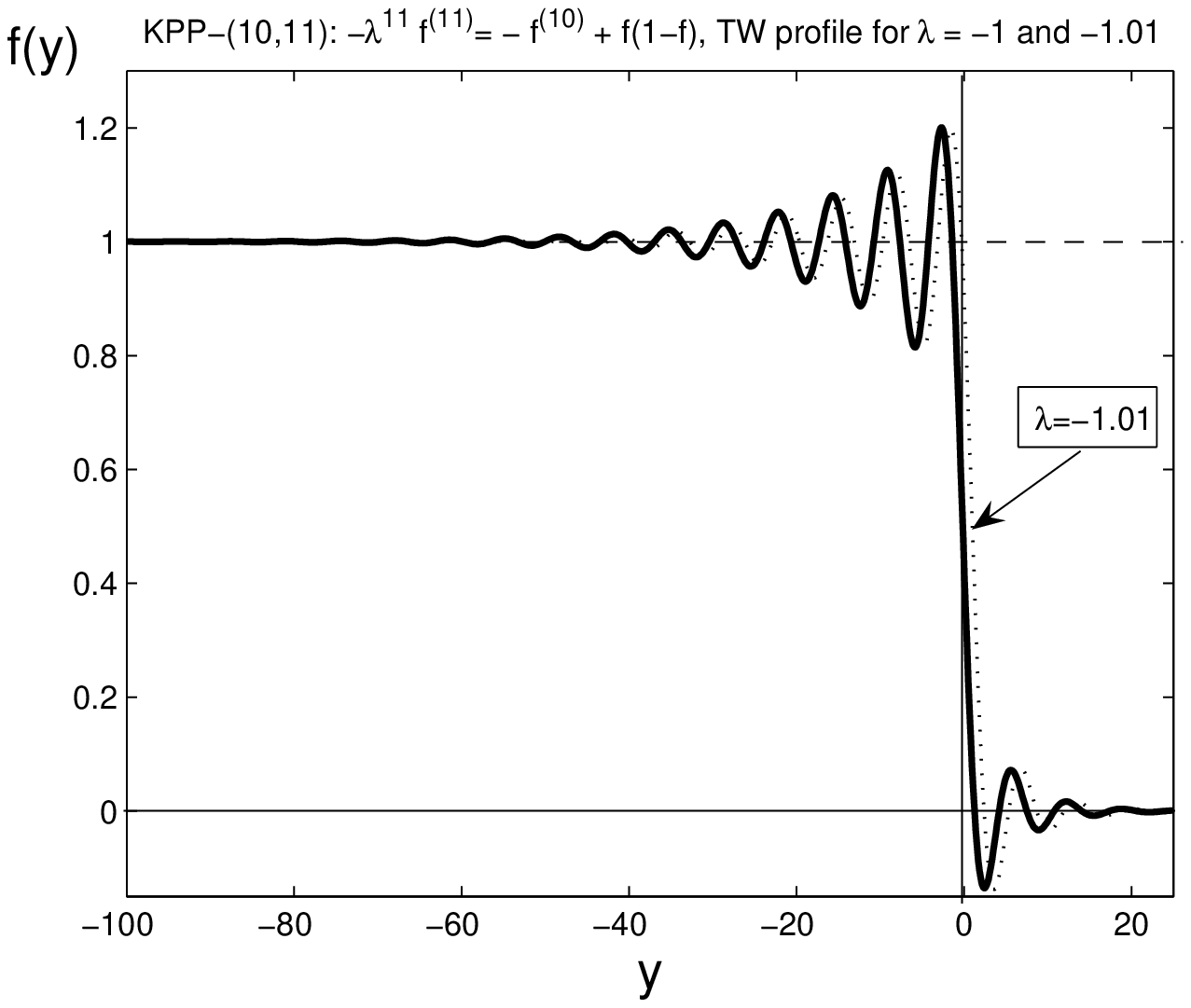}                        
}
 \vskip -.3cm
\caption{\rm\small Oscillatory one-sided ``periodic" TW patterns
for elliptic equations \ef{m32ell}--\ef{m36ell}.}
 \label{F77}
\end{figure}

In Figure \ref{F77}, we present TW profiles for the KPP--(10,11)
problem \ef{10.11} for $\l=0.2$ (a) (after a proper reflection,
this profile is very close to the stationary one in Figure
\ref{F37Period}) and for $\l=-1$, $-1.1$ (b).

\begin{figure}
\centering
\includegraphics[scale=0.70]{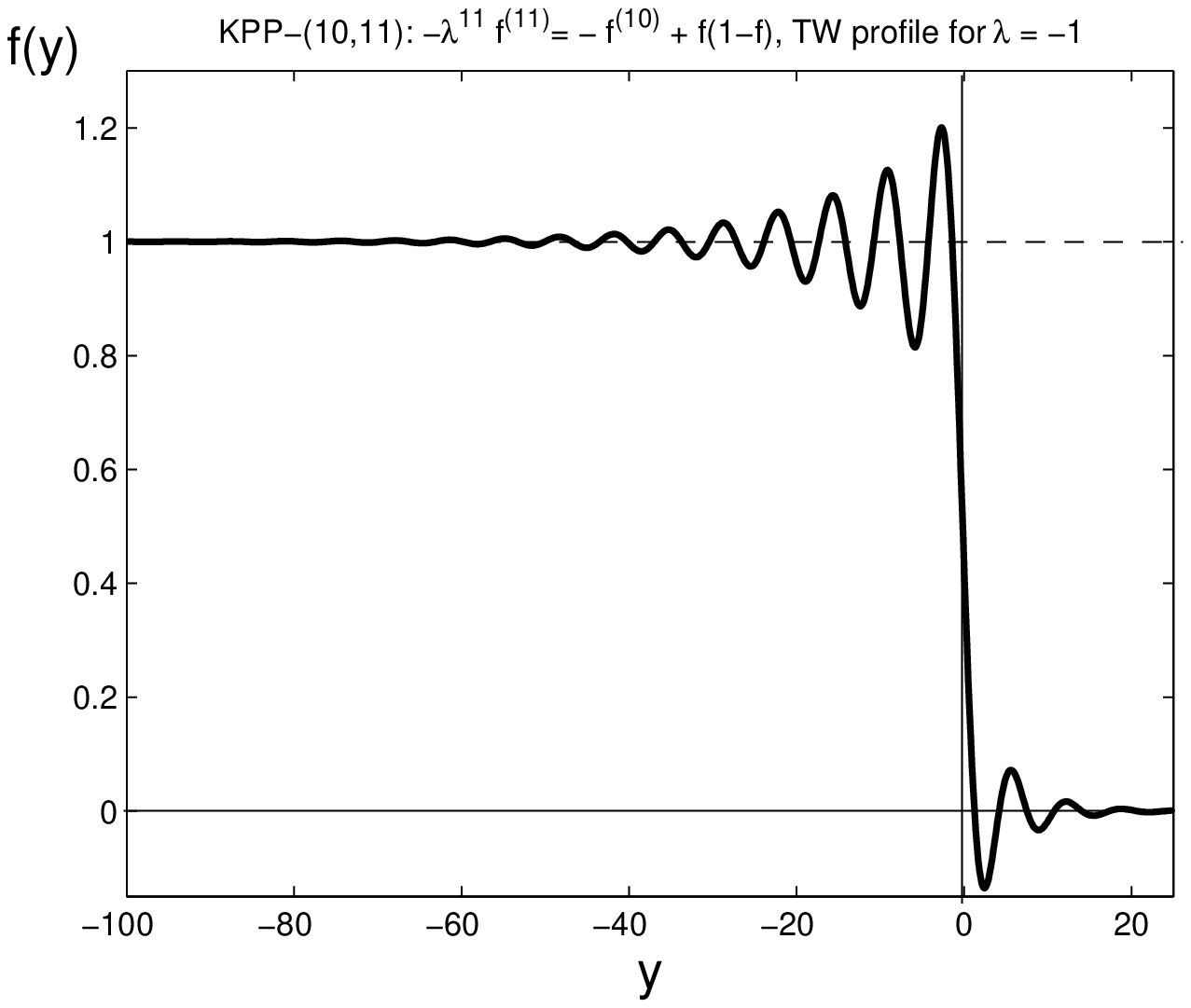}  
\vskip -.3cm
  \caption{A TW profile for the KPP--(11,12), \ef{11.12}, for $\l=1$ and 0.8.}
 \label{F88}
\end{figure}

In Figure \ref{F88}, we show  TW profiles for the PDE \ef{11.12}
for $\l=1$ and 0.8.


\section{Some semilinear PDEs that are eleventh-order in $t$}
 \label{S.1-4.11}

\begin{figure}
\centering \subfigure[(\ref{1.11}), $\l=-1,-2$]{
\includegraphics[scale=0.52]{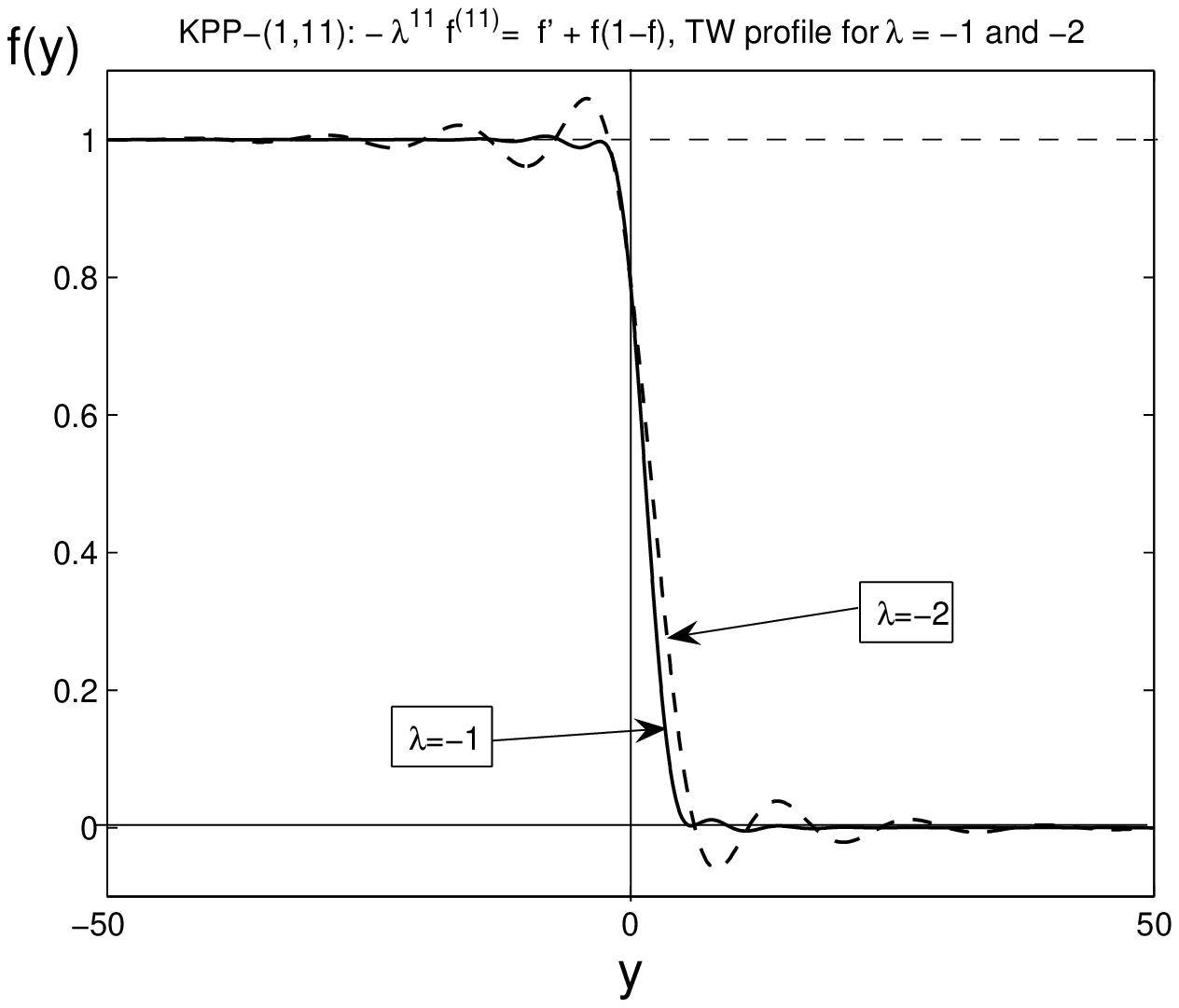}             
} \subfigure[(\ref{2.11}), $\l=-1$]{
\includegraphics[scale=0.52]{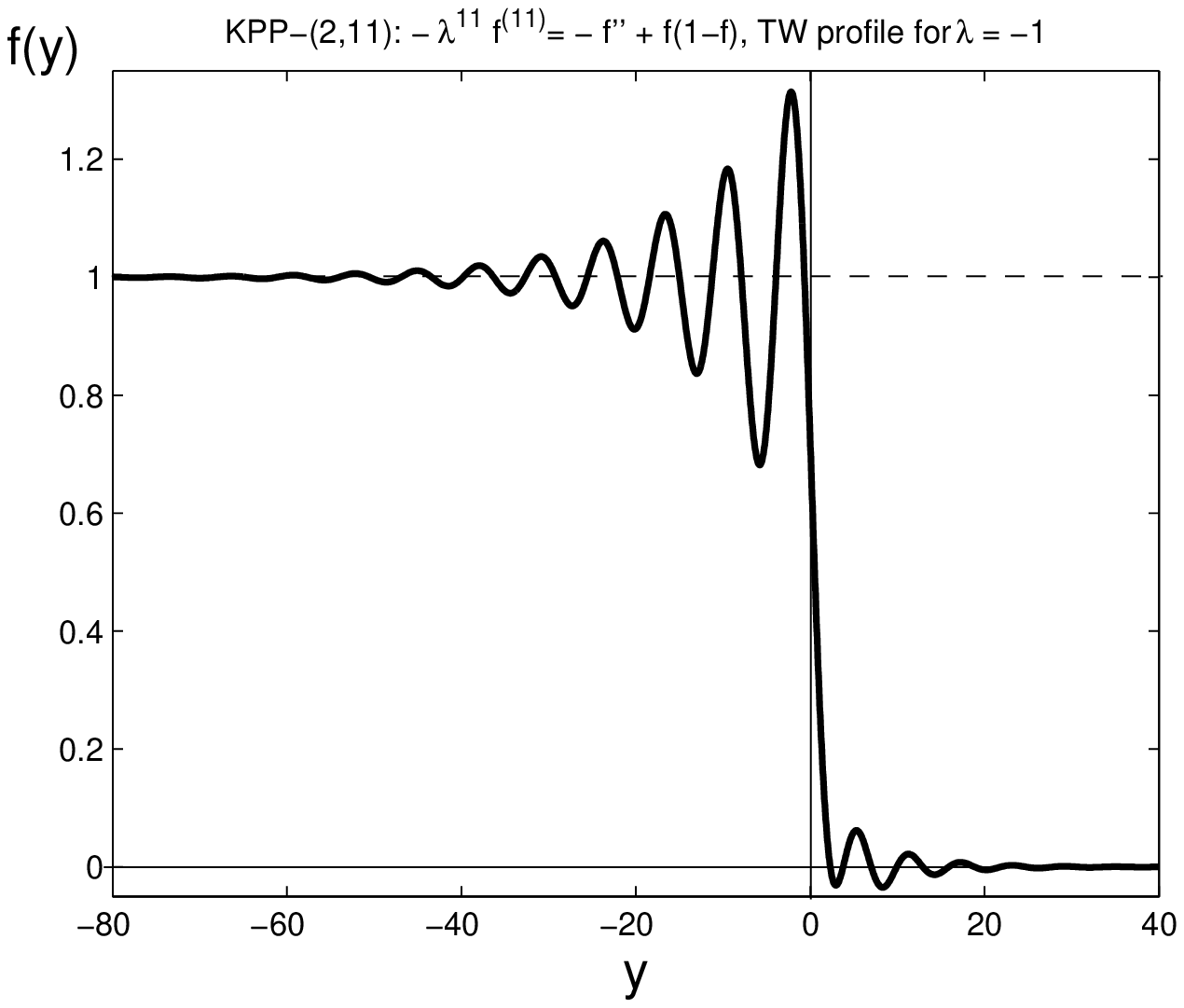}                        
} \subfigure[(\ref{3.11}), $\l=-1$]{
\includegraphics[scale=0.52]{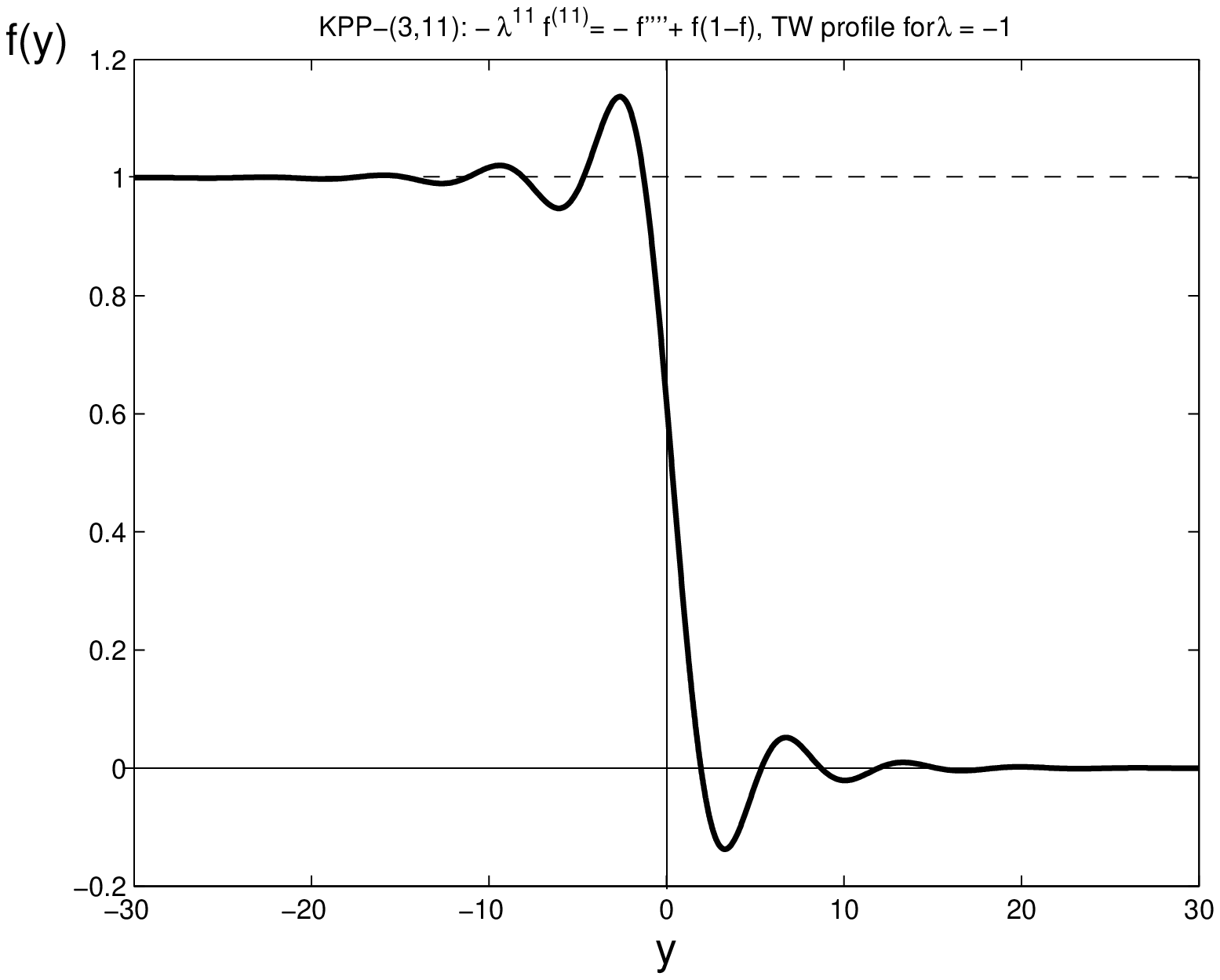}             
} \subfigure[(\ref{4.11}), $\l=-1$]{
\includegraphics[scale=0.52]{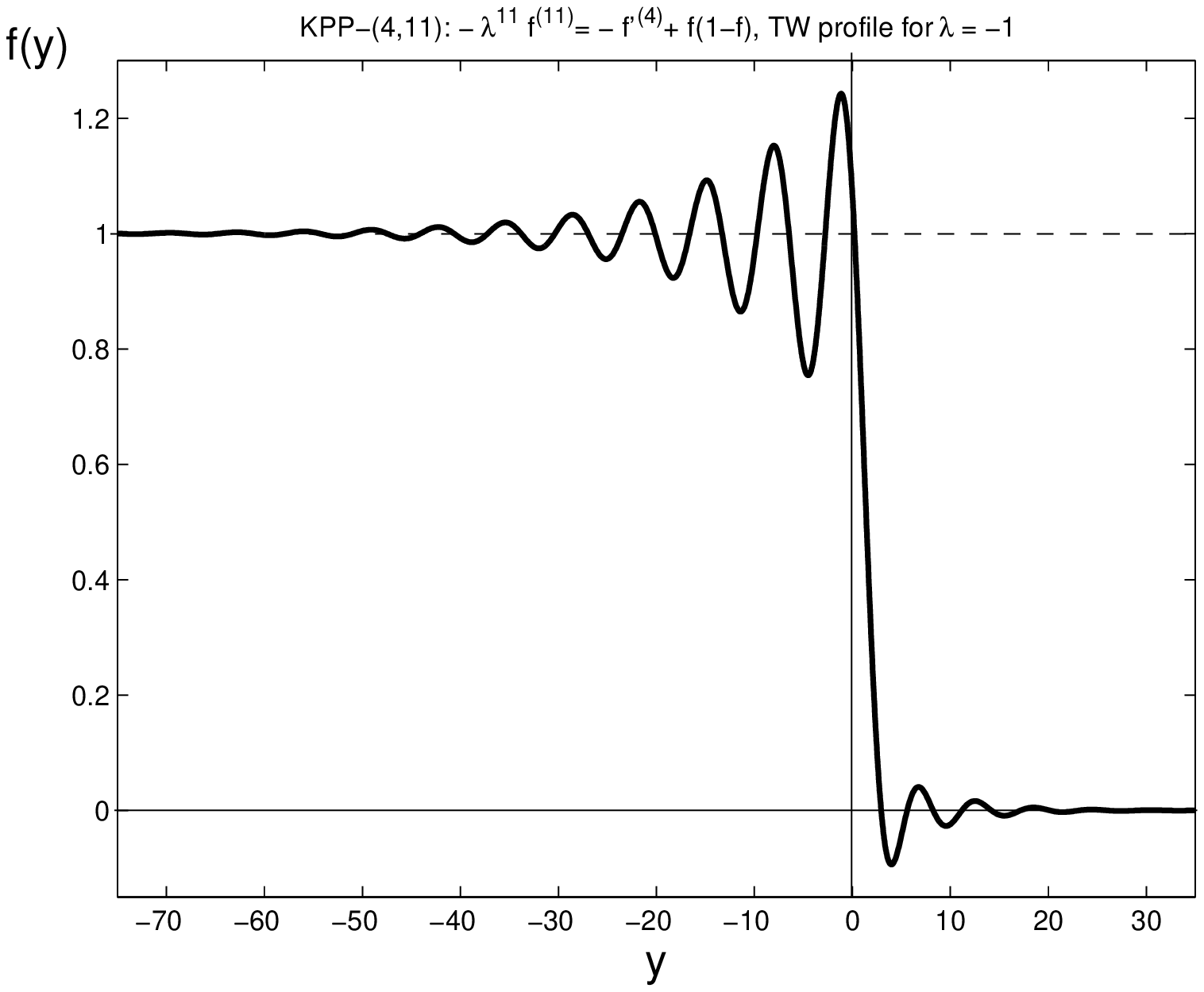}                        
}
 \vskip -.3cm
\caption{\rm\small  TW profiles for  equations
\ef{1.11}--\ef{4.11}.}
 \label{F.1-4.11}
 \end{figure}

For equations \ef{1.11}--\ef{4.11}, which all of different types,
the TW profiles look rather similar and are presented in Figure
\ref{F.1-4.11}. Note that, in all the cases, we fix negative
$\l=-1$. Choosing positive speeds $\l=+1$ always led to highly
oscillatory behaviour at the right-hand end point of the interval
of integration and no convergence was observed.

Note also that the first equation \ef{1.11} in this list admits
obvious (and well known) {\em explicit stationary solutions} for
$\l=0$:
 \be
 \label{1.11expl}
  \tex{
 (\ref{1.11}), \,\, \l=0: \quad f'=-f(1-f) \LongA f(y)=
 \frac{\eee^{-y}}{1+\eee^{-y}}.
 }
 \ee
 Therefore, one can expect a {\em branching} of TW profiles for
 small $|\l|>0$ from \ef{1.11} at $\l=0$. Then, if this is not a
 subcritical pitchfork $\l$-bifurcation, we expect existence
 $f(y)$, at least, for all small $\l>0$ (thought, not extensible
 to $\l=+1$, as mentioned above).

\section{Quasilinear dispersion KPP--(11,1) problem}
 \label{Sn11}

Consider the quasilinear ODE in \ef{m31n}, with $n>0$. First of
all, we justify that the TW profiles are assumed to have finite
interfaces at some finite point $y=y_0>0$; see extra details for
such a functional setting in \cite{GKPPII} (the singular
conditions as $y \to -\iy$ remain principally the same as in
\ef{BC1}). It is clear that, as $y \to y_0^-$, the source term
$f(1-f)$ can be neglected, so that the asymptotics of $f(y)$ near
this finite interface is described by both leading higher-order
terms, so one needs to study the following asymptotic ODE:
 \be
 \label{as66}
 -\l f'=-(|f|^nf)^{(11)} \LongA  \l f=(|f|^nf)^{(10)} \quad
 (\l>0).
  \ee
  Similar to the approach in \cite{GKPPII} for quasilinear
  parabolic equations, we claim that this asymptotic behaviour is
  given by
   \be
   \label{as67}
    \tex{
   f(y)=(y_0-y)^{\g} \var(s) \whereA s=\ln(y_0-y) \to
   -\iy \asA y \to y_0^-,  \quad \g=\frac{10}n,
   }
    \ee
  where the {\em oscillatory component} $\var(s)$ satisfies another complicated tenth-order
  ODE:
  \be
  \label{as68}
  P_{10}[|\var|^n \var]= \l \var \inB \re.
  \ee
Here, the linear polynomial operator $P_{10}[\phi]$ belongs to the
family $\{P_k[\phi], \, k \ge 0\}$ of operators that are
  constructed by the iteration
  \be
  \label{LH.7}
  P_{k+1}[\phi]= (P_k[\phi])' + (\g-k)P_k[\phi] \,\,\, \mbox{for}
  \,\,\, k=0,1,... \, , \,\,\,
  P_0[\phi]=\phi.
   \ee
In particular, this yields:
 \be
  \label{P75}
 P_3[\phi]= \phi''' + 3(\g-1) \phi'' + (3 \g^2 - 6 \g
+2) \phi'
  + \,\g(\g-1)(\g-2)\phi;
 \ee
 \be
 \label{P75.1}
  \begin{aligned}
 P_4[\phi]= &
 \phi^{(4)} + 2(2 \g-3) \phi'''+ (6 \g^2-18 \g+11) \phi''
    +2(2 \g^3-9 \g^2\\
   & + \,11 \g -3) \phi'
    +  \g(\g-1)(\g-2)(\g-3) \phi;
    \end{aligned}
    \ee
    \be
    \label{P75.2}
 \begin{aligned}
  P_5[\phi]= & \phi^{(5)}+5(\g -2)\phi^{(4)}
+5  (2\g^2-8\g+7)\phi'''\\
 &  +\, 5 (\g-2)(2\g^2-8\g +5)\phi''
 + (5\g^4-40\g^3+105 \g^2 \\
& -\, 100\g+24)\phi'
  + \, \g(\g-1)(\g-2)(\g-3)(\g-4)\phi, \quad
\mbox{etc.}
 \end{aligned}
 \ee
 The operator $P_{10}$ in \ef{as68} is too ambiguous to present it here.

 The next main point now is that the ODE \ef{as68} for $\var$ admits a {\em
 periodic} solution $\var_*(s)$, which, together with its
 stability set as $s \to -\iy$ (including the obvious translations in $s$),
  represents the actual asymptotic bundle of all admissible
  solutions of the form \ef{as67} near the interface.

  Examples of such periodic solutions $\var_*(s)$ for operators
  $P_6$ in \ef{as68} can be found in \cite[p.~192]{GSVR}; see also similar examples
   for $P_5$ in \cite{GBl6}. The
  higher-order case of $P_{9}$ (the parabolic one) is given
  in \cite[p.~143]{GSVR}, etc. The $P_{10}$ case is no much different,
  though, since $\var_*(s)$ can be more unstable as $s \to +\iy$
  (its unstability manifold $s \to -\iy$, i.e., approaching the
  interface, becomes even more dimensional), so numerics may get more
  difficult.

\ssk



Examples of full global solutions of the ODE in \ef{as66} are
presented in a number of figures below, where, for convenience and
by obvious reasons, we represent the function
 \be
 \label{F113}
  \tex{
 F(y)=|f(y)|^n f(y) \LongA
   F^{(11)}= \frac \l{n+1}\, |F|^{-\frac n{n+1}}
  F'+|F|^{-\frac n{n+1}}F \big(1-|F|^{-\frac n{n+1}}F\big).
 }
 \ee
Namely,  Figure \ref{F11n1} shows the TW profile for $\l=1$ and
$n=1$. The next Figure \ref{F11n2} shows TW profiles $F(y)$ again
for $n=1$ for $\l=0.5$ and $\l=1.196$. The last value turns out to
be close to the maximal value $\l_{\rm max}$ and we get the
estimate
 \be
 \label{maxn1}
n=1: \quad 1.196 \le \l_{\rm max}(1)<1.197.
 \ee
Recall that, for the semilinear case $n=0$, $\l_{\rm max}(0)$ is
slightly larger; see \ef{dispMax}.

\begin{figure}
\centering
\includegraphics[scale=0.70]{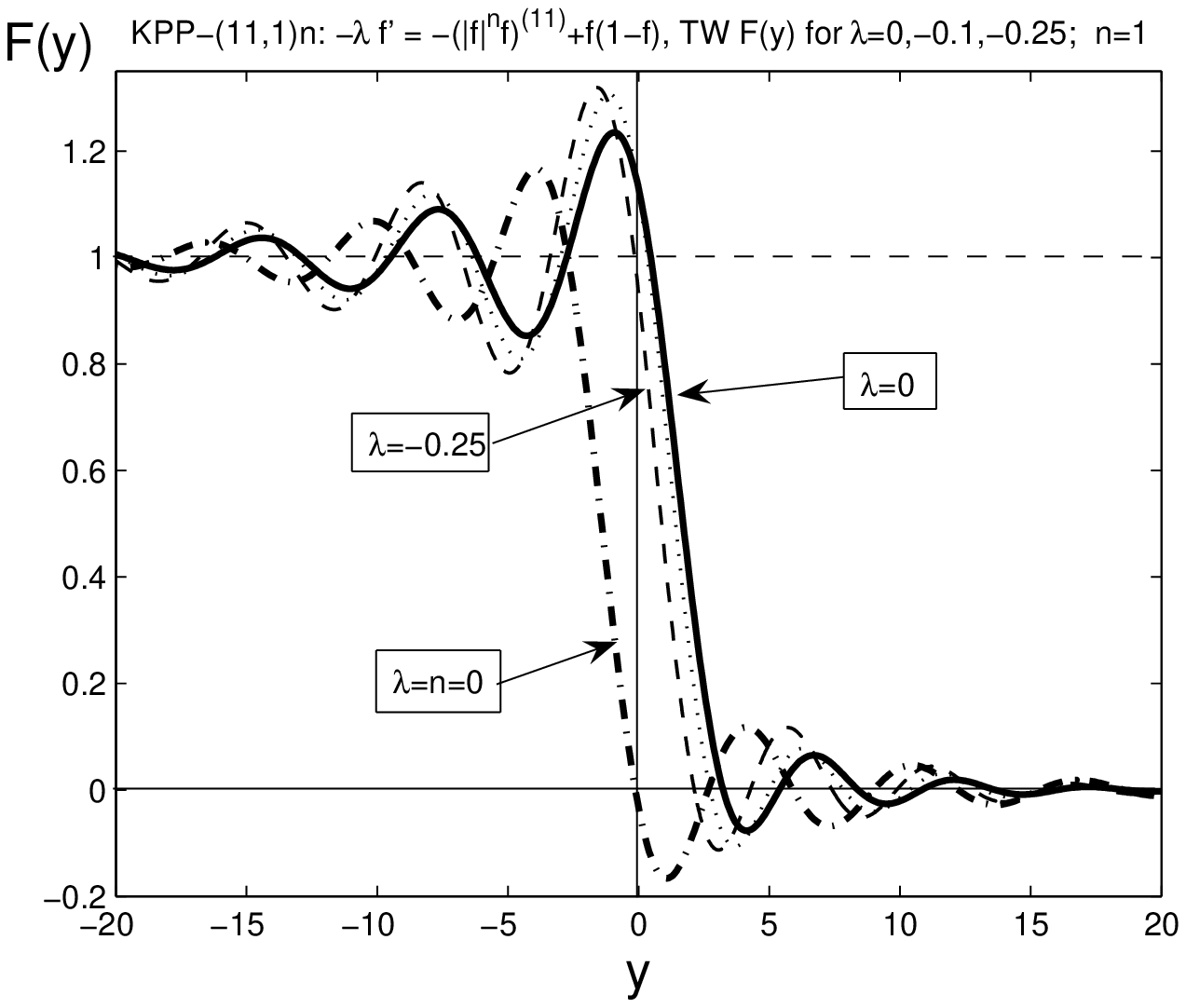}  
\vskip -.3cm
  \caption{TW profiles $F(y)$ satisfying \ef{F113}  for $\l=0, \,-0.1, \, -0.25$, and $n=1$.}
 \label{F11n3}
\end{figure}

TW profiles for some negative velocities $\l$ and $n=1$ are
presented in Figure \ref{F11n3}. For the sake of comparison, we
also indicate therein the stationary profile, with $\l=0$, for
$n=0$.

\begin{figure}
\centering
\includegraphics[scale=0.70]{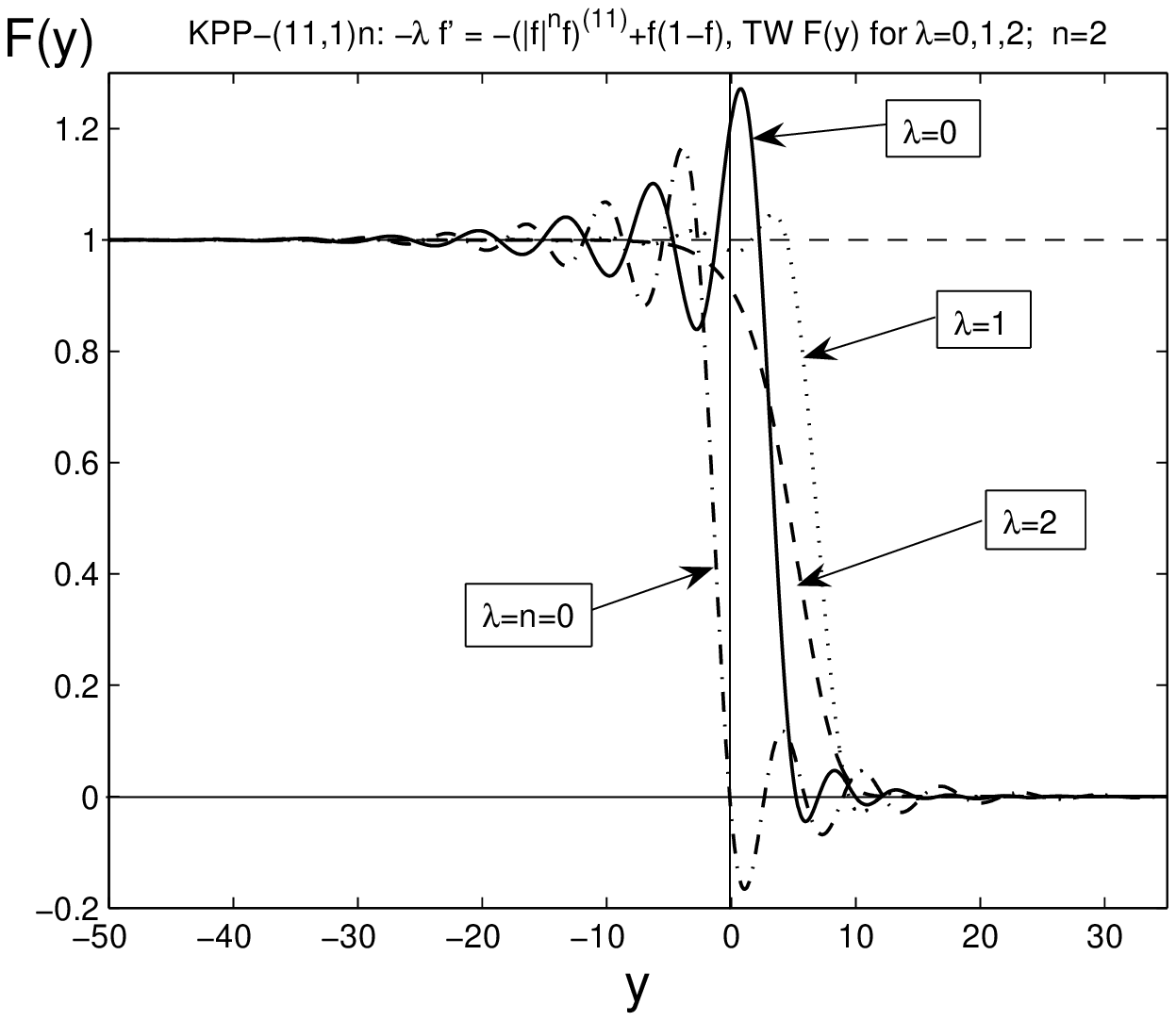}  
\vskip -.3cm
  \caption{TW profiles $F(y)$ satisfying \ef{F113}  for $\l=0,\, 1, \, 2$, and $n=2$.}
 \label{F11n4}
\end{figure}

\begin{figure}
\centering
\includegraphics[scale=0.70]{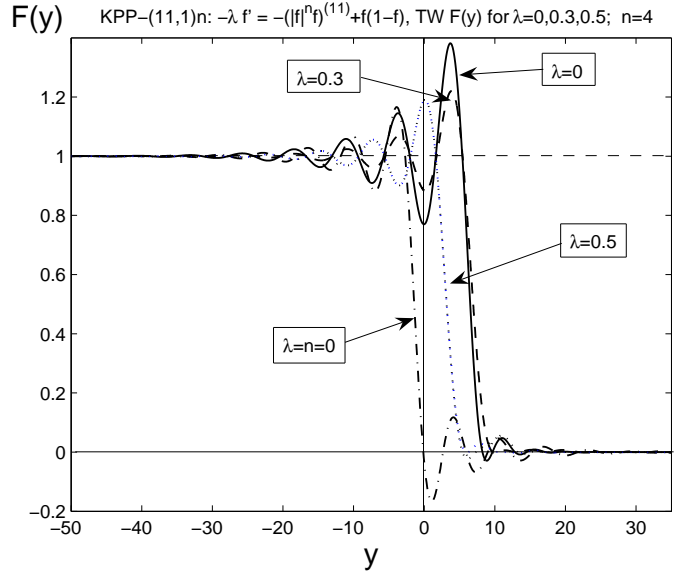}  
\vskip -.3cm
  \caption{TW profiles $F(y)$ satisfying \ef{F113}  for $\l=0$ (stationary),\, 0.3, \, 0.5, and $n=4$.}
 \label{F11n5}
\end{figure}

In Figure \ref{F11n4}, the TW profiles $F(y)$ correspond to a
larger $n=2$ and $\l=0,1,2$. The last value is not that far from
the maximal value: our computations show that
 \be
 \label{lmaxn2}
 n=2: \quad 2.25 \le \l_{\rm max}(2)< 2.26.
  \ee
For $n=4$, Figure \ref{F11n5} shows the profiles $F$ for $\l=0, \,
0.3$, and $0.5$.


 \begin{figure}
\centering
\includegraphics[scale=0.70]{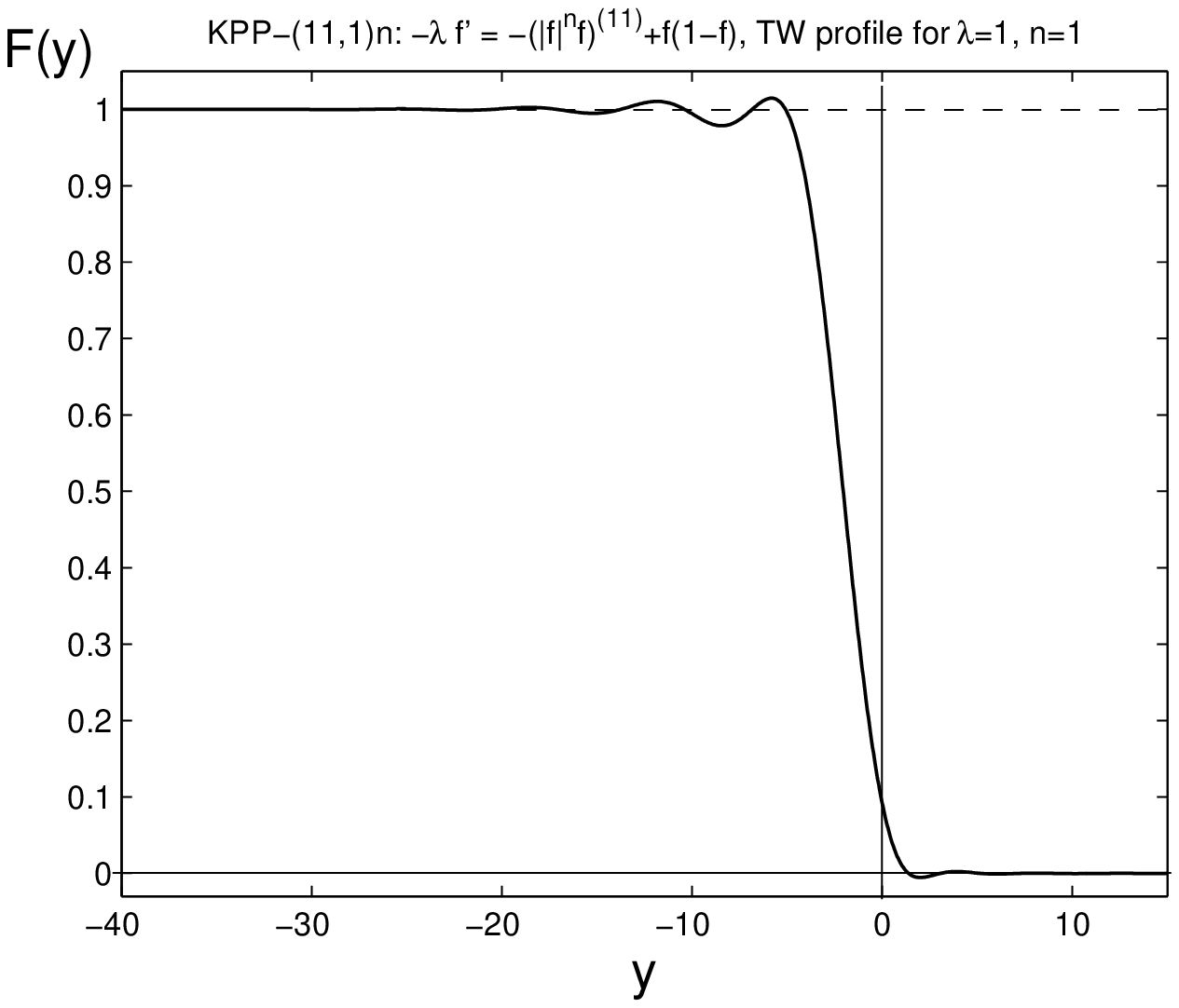}  
\vskip -.3cm
  \caption{A TW profile $F(y)$ satisfying \ef{F113}  for $\l=1$ and $n=1$.}
 \label{F11n1}
\end{figure}



 \begin{figure}
\centering
\includegraphics[scale=0.70]{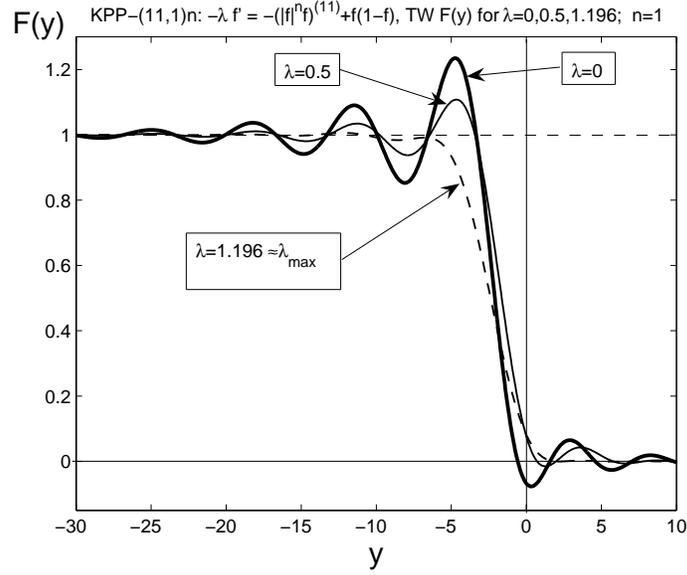}  
\vskip -.3cm
  \caption{A TW profile $F(y)$ satisfying \ef{F113}  for $\l=0, \,0.5, \, 1.196$, and $n=1$.}
 \label{F11n2}
\end{figure}



\section{A $\log t$-shift in the dispersion  KPP--(11,3) problem:
a centre
 subspace pattern}
 \label{Sdisplog}

Thus, we begin with a formal analysis of a kind of a ``centre
subspace behaviour", which generates a necessary $\log t$-shift of
the wave front. We restrict to a semilinear PDE. For similar
applications to  quasilinear (parabolic) ones, see \cite{GKPPII}.

\ssk

As a typical PDE, bearing  in mind  \ef{m33}, consider a
semilinear KPP-type problem for a PDE,
 \be
 \label{K1}
  u_{ttt}=\AAA u+u(1-u) \inB \re \times \re_+,
  \ee
  where $\AAA$ is a proper homogeneous isotropic and translational invariant linear differential
  operator satisfying some extra conditions specified below.
  For instance, we can fix the dispersion operator
   \be
   \label{disp11}
   \AAA=-D_x^{11} \quad (\mbox{cf. (\ref{m31})}).
   \ee
   We assume that the corresponding ODE problem
   \be
   \label{K2}
   -\l_0^3 f'''= \AAA f+f(1-f),
   \ee
   with the conditions \ef{BC1} admits a unique solution $f$.

   Attaching the solution $u(x,t)$ to the
    front moving and setting, as usual,  $x_f(t) \equiv \l_0t -g(t)$,
    the PDE reads
    \be
    \label{K21}
    u(x,t)=v(y,t), \quad y=x- \l_0 t+g(t).
     \ee
 Then $v$ satisfies the following perturbed equation:
  \be
  \label{eq55}
   \begin{aligned}
    &v_{ttt} + 3 v_{tty}(-\l_0+g')' + 3v_{tyy}(-\l_0+g')^2 +
 3 v_{ty}g'' \ssk \\
 = &\AAA v + v(1-v)  - v_{yyy}(-\l_0+g')^3 -
3 v_{yy}(-\l_0+g')^2 g'' - v_y g'''.
 \end{aligned}
 \ee

 Thus, again, as usual, we assume that $g'(t) \to 0$ as
 $t \to +\iy$ sufficiently fast, i.e., at least algebraically, so
 that
  \be
  \label{gg1}
  |g''(t)| \ll |g'(t)| \andA |g'''(t)| \ll |g''(t)|, \,\,\, \mbox{etc.} \forA t \gg 1.
  \ee

   We next ``linearize" \ef{eq55} by setting
    \be
    \label{K3}
    v(y,t)=f(y)+w(y,t).
     \ee
     Then, using \ef{K2},
      yields the following perturbed equation, where,
     according to \ef{gg1}, we keep the leading terms only:
      \be
      \label{gg2}
       \begin{aligned}
       w_{ttt}- & 3\l_0 w_{tty}+3\l_0^2 w_{tyy}  + 3 w_{ty}g''
      = \BB w - 3 \l_0^2f'' g''-f'g'''-w^2, \ssk \\
      & \mbox{where} \quad
      \BB w= \AAA w + (1-2f)w
       + \l_0^3 w_{yyy}.
       \end{aligned}
       \ee
   Note that here we face an essentially non-autonomously perturbed flow, so we cannot use advanced semigroup
   theory; cf. \cite{Lun}. Instead, we will apply  formal asymptotic
   expansion techniques.

Now, assuming that, in this $g(t)$-moving frame, there exists the
 convergence as in \ef{TW3}, so that $w(t) \to 0$ as $t \to +\iy$,
 one can see that, under the hypothesis \ef{gg1}, the leading non-autonomous
 perturbations
in \ef{gg2} are those of order $O(g''(t))$, since the rest of the
terms are negligible as $t \to +\iy$. Therefore, one needs to
balance these major terms, but then the actual behaviour of $g(t)$
for $t \gg1$ (and, hence, proper $\log t$-shifts of the front)
will depend on the next matching.

Thus, under the hypothesis \ef{gg1}, the only possible way to
balance {\em all} the terms therein (including the quadratic one
$-w^2$) for $t \gg 1$ is to assume the following expansion:
 \be
 \label{gg3}
 w(y,t) = g'(t) \psi(y)+ \e(t) \var(y)...\asA t \to +\iy
 \whereA |\e(t)| \ll |g'(t)|
  \ee
  is still an unknown coefficient.
Substituting \ef{gg3} into \ef{gg2} yields
 \be
 \label{gg31}
  \begin{aligned}
 &g^{(4)}\psi + \e'''\var +3 (g''' \psi'+\e'' \var')(-\l_0+g') \\
 &+
 3(g'' \psi''+\e'\var'')(\l_0^2-2\l_0g'+(g')^2)+3(g'' \psi'+\e'\var')g''+...\\
 &
 = g'\BB \psi +\e \BB \var- (g')^2 \psi^2- 2\e g' \var \psi\\
 & + \e^2
 \var^2-(f'''+g'\psi'''+\e \var''')(-\l_0^3+3\l_0^2g'-3
 \l_0(g')^2+(g')^3)
 \\
 &-3(f''+g' \psi''+\e \var'')(-\l_0+g')g''-(f'+
 g' \psi'+\e \var')g''+... \,,
 \end{aligned}
  \ee
  where we have omitted some obviously negligible terms.

Using \ef{gg1} and \ef{gg3} in
 balancing first the leading terms of the order
  $O(g'(t))$ yields the inhomogeneous equation for $\psi$:
   \be
   \label{gg4}
    \tex{
 O(g'(t))\,\,\big(=O(\frac 1t), \,\,\mbox{see below}\big): \quad  \BB \psi - 3\l_0^2 f'''=0.
   }
   \ee
   Then balancing the rest of the terms in \ef{gg31} requires
 \be
 \label{hh41}
  \tex{
 g''(t) \sim -(g'(t))^2 \sim \e(t), \,\,\, \mbox{i.e.,} \,\,\,
  g(t)= k \log t, \,\, \mbox{$g'(t)= \frac kt, \,\,
 g''(t)=- \frac k{t^2}$,} \,\, \e(t)= \frac 1{t^2}.
 }
 \ee
 Then, we obtain the second inhomogeneous singular
 Sturm--Liouville problem for $\var$:
  \be
  \label{SL1}
   \tex{
  O\big(\frac 1{t^2}\big): \quad
  \BB \var= k(3 \l_0 f''-3 \l_0^2 \psi''-f') + k^2(\psi^2+3
  \l_0^2\psi'''-3\l_0 f''').
 }
   \ee

    Thus, the first simple asymptotic ODE in \ef{hh41} gives the $\log
    t$-dependence as in \ef{1.3}. Finally, we arrive at the
    following system for $\{\psi,\var\}$:
     \be
     \label{sys21}
     \left\{
     \begin{aligned}
& \BB \psi = 3\l_0^2 f'', \\
   &
   \BB \var=  k(3 \l_0 f''-3 \l_0^2 \psi''-f') + k^2(\psi^2+3
  \l_0^2\psi'''-3\l_0 f''').
    \end{aligned}
     \right.
 \ee
  Solving this system, with typical boundary conditions as in
  \ef{BC1}, allows then continue the expansion of the solutions of
  \ef{gg2}
  close to an ``affine centre subspace" of $\BB$ governed by the obvious (by translation)
   spectral pair
   \be
   \label{pair1}
   \hat \l_0=0 \andA \hat \psi_0(y)= f'(y).
    \ee
 The asymptotic expansion for $t \gg 1 $ then takes the form
  \be
  \label{as22}
  \tex{
  w(y,t)= \frac kt \, \psi(y) + \frac 1{t^2} \,\var(y) + ...\, ,
  }
  \ee
  which can be easily extended by introducing further terms, with
  similar inhomogeneous Sturm--Liouville problems for the
  expansion coefficients.

Since $\BB$ does not have a discrete spectrum, we cannot derive a
simple algebraic equation for $k$ by demanding the standard
orthogonality of the right-hand side in the second equation in
\ef{sys21} to the adjoint eigenvector $\hat \psi^*_0$ of $\BB^*$
in the $L^2$-metric (in which the adjoint operator $\BB^*$ is
obtained), like
 \be
 \label{as23}
k: \quad   \langle k(3 \l_0 f''-3 \l_0^2 \psi''-f') + k^2(\psi^2+3
  \l_0^2\psi'''-3\l_0 f'''), \, \hat \psi_0^* \rangle =0.
 \ee
 Therefore,  the system \ef{sys21} cannot itself  determine
    the actual value of $k$ therein. As we have mentioned, the latter requires a
    difficult matching analysis in Inner and Outer Regions, which,
    for all present KPP--problems, remains an open
    problem.

    \ssk

    One can see that the above elementary conclusion well
    corresponds to a ``centre subspace analysis" of the
    non-autonomous PDE \ef{gg2}, and then $\t = \log t$ naturally
    becomes the corresponding ``slow" time variable; see various
    examples in \cite{AMGV} of such slow motion along centre subspaces in
    nonlinear parabolic problems with global and blow-up
    solutions.
    In the latter case, the slow time variable is $\t = - \ln
    (T-t) \to +\iy$ as $t$ approaches finite blow-up time $T^-$.

\ssk

 For the semilinear higher-order reaction-absorption equations
such as
 \be
 \label{utp}
 u_t=-u_{xxxx}-|u|^{p-1}u \inB \re \times \re_+ \withA p>1,
 \ee
existence of $\log t$-perturbed global asymptotics was established
in \cite{GalCr}. For finite-time extinction, with $-1<p<1$ in
\ef{utp}, this was done in \cite{Galp1}. For the corresponding
blow-up problem with the combustion source
 $$
 u_t=-u_{xxxx}+|u|^{p-1}u \inB \re \times \re_+ \withA p>1,
 $$
$\log(T-t)$-dependent blow-up singularities were constructed in
\cite{Gal2m1}. We must admit that any justification of such $\log
t$-front corrections (and finding corresponding classes of initial
data) in any of KPP--$(k,l)$ problems with $k,l
> 2$ is much more difficult and remains an open problem.

\section{Final conclusions}

{\bf 1.} In all of KPP--$(k,l)$
 problems, there exist TW solutions
for various values of the speeds $\l$. Most of them satisfy
singular boundary conditions at $\pm \iy$. However, for some types
of PDEs involved, these can be oscillatory and/or periodic either
as $y \to +\iy$, or $y \to -\iy$, and hence require special
setting.

\ssk

{\bf 2.} In all of the higher-order KPP--$(k,l)$ problems, we did
not observe the classic KPP--2 phenomenon for \ef{1.1}, \ef{TW1}
\cite{KPP} of existence of the {\em minimal} speed $\l_0=2$, such
that TWs exist for $\l \ge 2$ only. It seems that this phenomenon
is directly connected with the Maximum Principle (and other
features of order-preserving flows), and becomes non-generic and
non-existent when this fails.

\ssk

 {\bf 3.} Moreover, in several KPP--$(k,l)$ problems, on the
contrary, we observed existence of a {\em maximal} speed $\l_{\rm
max}$, so that for slightly $\l \ge \l_{\rm max}$, TW profiles do
not exist. In particular, this always happens for parabolic
problems, \cite{GKPPI}.

\ssk

{\bf 4.} In all the KPP--$(k,l)$ problems, there exists a formal
justification of existence of $\log t$-drift of the propagating
fronts  in PDE setting along an ``(affine) centre subspace" of
semilinear rescaled operators involved. Then $\t=\log t$ naturally
appears as the corresponding ``slow" time variable.

\ssk

{\bf 5.} We must admit that the important PDE problem on the
actual structure of the omega-limit set $\o(H)$ of the rescaled
(properly shifter) solution $u(x,t)$ of the various KPP-problems
with the Heaviside data $H(-x)$ remains open. In particular, it is
not still known whether $\o(H)=\{f(\cdot;\l_0)\}$ for some $\l_0
\in \Lambda$, i.e., whether $\o(H)$ consists of a single TW
profile.
 We believe that this
problem deserves further study by analytical and PDE numerical
methods, but expect it to be very difficult.


\end{document}